\documentclass[twoside,reqno]{amsart}

\usepackage{amsmath,amssymb,amsbsy,amstext,amsxtra,latexsym}
\usepackage{graphicx}
\usepackage{subfig}
\usepackage{enumerate}
\usepackage[all]{xy}
\usepackage{url}

\usepackage{hyperref}

\hypersetup{%
pdftitle={A complex surface of general type with $p_g=0$, $K^2=2$ and
 $H_1=\mathbb{Z}/4\mathbb{Z}$},%
pdfauthor={Heesang Park, Jongil Park, Dongsoo Shin},%
pdfkeywords={$\mathbb{Q}$-Gorenstein smoothing, rational blow-down surgery,
 surface of general type},%
}

\newtheorem{theorem}{Theorem}[section]

\newtheorem*{maintheorem}{Theorem}
\newtheorem{lemma}[theorem]{Lemma}
\newtheorem{proposition}[theorem]{Proposition}

\theoremstyle{definition}

\theoremstyle{remark}

\newtheorem{remark}[theorem]{Remark}

\numberwithin{equation}{section}

\newcommand{\uc}[1]{\ensuremath \overset{#1}{\circ}}
\newcommand{\udc}[2]{\ensuremath \underset{#2}{\overset{#1}{\circ}}}

\DeclareMathOperator{\im}{Im}

\DeclareMathOperator{\divisor}{Div}

\def\sheaf#1{\ensuremath \mathcal#1}

\begin{document}

\title[A surface with $p_g=0$, $K^2=2$ and $H_1=\mathbb{Z}/4\mathbb{Z}$]{A complex surface of general type \\ with $p_g=0$, $K^2=2$ and $H_1=\mathbb{Z}/4\mathbb{Z}$}

\author{Heesang Park}

\address{School of Mathematics, Korea Institute for Advanced Study, Seoul 130-722, Korea}

\email{hspark@kias.re.kr}

\author{Jongil Park}

\address{Department of Mathematical Sciences, Seoul National University, Seoul 151-747, Korea \& Korea Institute for Advanced Study, Seoul 130-722, Korea}

\email{jipark@snu.ac.kr}

\author{Dongsoo Shin}

\address{Department of Mathematics, Chungnam National University, Daejeon 305-764, Korea}

\email{dsshin@cnu.ac.kr}


\subjclass[2000]{Primary 14J29; Secondary 14J10, 14J17, 53D05}

\keywords{$\mathbb{Q}$-Gorenstein smoothing, rational blow-down surgery, surface of general type}

\begin{abstract}
We construct a new minimal complex surface of general type with $p_g=0$, $K^2=2$ and $H_1=\mathbb{Z}/4\mathbb{Z}$ (in fact $\pi_1^{\text{alg}}=\mathbb{Z}/4\mathbb{Z}$),
which settles the existence question for numerical Campedelli surfaces
with all possible algebraic fundamental groups.
The main techniques involved in the construction are a rational blow-down surgery and a $\mathbb{Q}$-Gorenstein smoothing theory.
\end{abstract}

\maketitle

\section{Introduction}

One of the fundamental problems in the classification of complex surfaces is to find a new family of complex surfaces of general type with $p_g=0$. In this paper we construct a new minimal complex surface of general type with $p_g=0$, $K^2=2$ and $H_1=\mathbb{Z}/4\mathbb{Z}$.

In order to classify complex surfaces of general type with $p_g=0$, it seems to be natural to classify them first up to their topological types. For instance if $S$ is a \emph{numerical Godeaux surface}, that is, a minimal complex surface of general type with $p_g=0$ and $K^2=1$, then it is known that $\pi_1^{\text{alg}}(S) = \mathbb{Z}/m\mathbb{Z}$  for some $1 \le m \le 5$, where $\pi_1^{\text{alg}}(S)$ is the algebraic fundamental group of $S$; cf.\thinspace Reid~\cite{Reid-78}. It is conjectured that the moduli space of numerical Godeaux surfaces has exactly five irreducible components corresponding to each $\pi_1=\mathbb{Z}/m\mathbb{Z}$ for all $1 \le m \le 5$. The conjecture is proved for $m \ge 3$; cf.\thinspace Reid~\cite{Reid-78}. Furthermore for each $1 \le m \le 5$ there are numerical Godeaux surfaces with $\pi_1=\mathbb{Z}/m\mathbb{Z}$; cf.\thinspace Bauer-Catanese-Pignatelli~\cite{Bauer-Catanese-Pignatelli}.

In case of \emph{numerical Campedelli surfaces}, that is, minimal complex surfaces of general type with $p_g=0$ and $K^2=2$, it has been known by Reid~\cite{Reid} and Xiao~\cite{Xiao-85} that the algebraic fundamental group of a numerical Campedelli surface is a finite group of order $\le 9$. Furthermore the topological fundamental groups $\pi_1$ for any numerical Campedelli surfaces are also of order $\le 9$ in as far as they have been determined. Hence it is a natural conjecture that $\lvert \pi_1 \rvert \le 9$ for all numerical Campedelli surfaces.

Conversely one may ask whether every group of order $\le 9$ occurs as the topological fundamental group or as the algebraic fundamental group of a numerical Campedelli surface. It is proved that the dihedral groups $D_3$ of order $6$ or $D_4$ of order $8$  cannot be fundamental groups of numerical Campedelli surfaces; Naie~\cite{Naie}, Mendes Lopes-Pardini~\cite{Mlopes-Pardini}, Mendes Lopes-Pardini-Reid~\cite{Mlopes-Pardini-Reid}, Reid~\cite{Reid}.

However all groups of order $\le 9$, \emph{except} $D_3$, $D_4$, $\mathbb{Z}/4\mathbb{Z}$, $\mathbb{Z}/6\mathbb{Z}$, occur as the topological fundamental groups of numerical Campedelli surfaces. The first example of numerical Campedelli surfaces was constructed by Campedelli himself and his example has $\pi_1=(\mathbb{Z}/2\mathbb{Z})^{\oplus 3}$. Furthermore Mendes Lopes, Pardini, and Reid (\cite{Mlopes-Pardini}, \cite{Mlopes-Pardini-Reid}, \cite{Reid}) constructed and classified all numerical Campedelli surfaces with $\lvert \pi_1^{\text{alg}} \rvert=8, 9$ and showed that the topological fundamental group equals the algebraic fundamental group. An example with $\pi_1=\mathbb{Z}/7\mathbb{Z}$ was constructed by Reid~\cite{Reid-91}. Catanese~\cite{Catanese} constructed an example with $\pi_1=\mathbb{Z}/5\mathbb{Z}$ by taking a $\mathbb{Z}/5\mathbb{Z}$-quotient of a certain double cover of a linearly symmetric quintic. An example with $\pi_1=\mathbb{Z}/2\mathbb{Z} \oplus \mathbb{Z}/2\mathbb{Z}$ was constructed by Inoue~\cite{Inoue} as the quotients of the hypersurfaces of the product of three elliptic curves. Such an example was also constructed by J. Keum~\cite{Keum}. Recently, several examples of numerically Campedelli surfaces with $\lvert \pi_1 \rvert \ge 3$ were constructed by the so-called product-quotient method, that is, first take a product of two curves and then take a quotient of the product by a group action;  Bauer-Catanese-Grunewald-Pignatelli~\cite{Bauer-Catanese-Grunewald-Pignatelli}, Bauer-Pignatelli~\cite{Bauer-Pignatelli}. Classically many examples of numerical Campedelli surfaces with $\lvert \pi_1 \rvert \ge 3$ are constructed by the method of taking quotients of group action.
On the other hand numerical Campedelli surfaces with small fundamental groups were constructed by the so-called $\mathbb{Q}$-Gorenstein smoothing method developed in Y. Lee-J. Park~\cite{Lee-Park-K^2=2}. For instance a simply connected numerical Campedelli surface was constructed by Y. Lee-J. Park~\cite{Lee-Park-K^2=2} and an example with $\pi_1=\mathbb{Z}/2\mathbb{Z}$ was constructed by J. Keum-Y. Lee-H. Park~\cite{Keum-Lee-Park}.

Unlike the case of topological fundamental group, there is a numerical Campedelli surface with $H_1=\mathbb{Z}/6\mathbb{Z}$. Recently Neves and Papadakis~\cite{Neves-Papadakis} constructed a numerical Campedelli surface with $H_1=\mathbb{Z}/6\mathbb{Z}$
(in fact $\pi_1^{\text{alg}}=\mathbb{Z}/6\mathbb{Z}$)
by constructing the canonical ring of the  \'etale six to one cover using serial unprojection together with a suitable base point free action of $\mathbb{Z}/6\mathbb{Z}$. Therefore all abelian groups of order $\le 9$ \emph{except} $\mathbb{Z}/4\mathbb{Z}$ occur as the first homology groups (and algebraic fundamental groups) of numerical Campedelli surfaces. Nevertheless, the question on the existence of numerical Campedelli surfaces with a given topological type was completely open for $\mathbb{Z}/4\mathbb{Z}$ until now.

The main result of this paper is the following.

\begin{maintheorem}
\label{theorem:main-theorem-intro}
There exists a minimal complex surface of general type with $p_g=0$, $K^2=2$ and $H_1=\mathbb{Z}/4\mathbb{Z}$.
Furthermore its algebraic fundamental group $\pi_1^{\text{alg}}$ is also $\mathbb{Z}/4\mathbb{Z}$.
\end{maintheorem}

To our knowledge there were no such examples previously known, and it settles the existence question for numerical Campedelli surfaces with $\lvert H_1 \rvert \le 9$ and furthermore the existence question for $\lvert \pi_1^{\text{alg}} \rvert \le 9$ because it can be shown that the algebraic fundamental group of the numerical Campedelli surface in the main theorem is also equal to $\mathbb{Z}/4\mathbb{Z}$. The reason is the following argument used in Bauer-Catanese-Pignatelli~\cite{Bauer-Catanese-Pignatelli} (Pignatelli~\cite{Pignatelli}): As mentioned above, by Reid~\cite{Reid} and Xiao~\cite{Xiao-85}, any numerical Campedelli surface has $\lvert \pi_1^{\text{alg}} \rvert \le 9$. Then it follows that the algebraic fundamental group is the quotient of the topological fundamental group by the intersection $N$ of all normal subgroups of finite index. Since the first homology group $H_1$ of the main example is finite, the intersection $N$ is contained in the commutator subgroup of the topological fundamental group. It follows that $H_1$, the abelian quotient of the topological fundamental group, is also the abelian quotient of the algebraic fundamental group. It has been known by Reid~\cite{Reid} and Xiao~\cite{Xiao-85} again that the only possible non-abelian group of algebraic fundamental groups for a numerical Campedelli surface is the quaternion group $Q_8$; but, its abelian quotient is $\mathbb{Z}/2\mathbb{Z} \oplus \mathbb{Z}/2\mathbb{Z}$. Therefore the algebraic fundamental group of the surface in the main theorem above is abelian and hence it is equal to $\mathbb{Z}/4\mathbb{Z}$.

In order to construct such a surface we use the $\mathbb{Q}$-Gorenstein smoothing method in Y. Lee-J. Park~\cite{Lee-Park-K^2=2}. We blow up an elliptic Enriques surface $Y$ in a suitable set of points so that we obtain a surface $Z$ with special linear chains of $\mathbb{P}^1$'s.
Then, by contracting these linear chains, we obtain a singular surface $X$ with permissible singular points in Section~\ref{section:main-construction}.
Each of these singular points admits a local $\mathbb{Q}$-Gorenstein smoothing. In Section~\ref{section:global-Q} we show that the obstruction space of a global $\mathbb{Q}$-Gorenstein smoothing of the singular surface $X$ is zero by a similar strategy as in J. Keum-Y. Lee-H. Park~\cite{Keum-Lee-Park}. Hence these local smoothings can be glued to a global $\mathbb{Q}$-Gorenstein smoothing of the whole singular surface $X$. Finally, we show that a general fiber $X_t$ of a $\mathbb{Q}$-Gorenstein smoothing of $X$ is the desired surface.

The main ingredient of this paper is the proof of $H_1(X_t)=\mathbb{Z}/4\mathbb{Z}$ in Section~\ref{section:Z_4}. For  proving it, instead of a general fiber $X_t$, we consider a rational blow-down $4$-manifold $\overline{Z}$ obtained by replacing certain small neighborhoods of the linear chains in $Z$ with the corresponding Milnor fibers. The rational blow-down $4$-manifold $\overline{Z}$ is known to be diffeomorphic to a general fiber $X_t$ by Milnor fiber theory. We first prove $\lvert H_1(\overline{Z};\mathbb{Z}) \rvert = 4$ by analyzing the long exact sequence of homologies induced from a pair ($\overline{Z}$, the Milnor fibers). We then construct carefully a certain loop $\beta$ lying on a certain curve of genus $2$ in the blown-up surface $Z$ and we prove that $H_1(\overline{Z};\mathbb{Z})$ is generated by the loop $\beta$, which finishes the proof. In order to show that the loop $\beta$ is a generator of $H_1(\overline{Z};\mathbb{Z})$, we lift the loop $2\beta$, a loop twice of $\beta$, up to a certain blown-up K3 surface $W$ which is a unramified double cover of the blown-up Enriques surface $Z$ and in $W$ we construct a real two dimensional surface $U_T$ such that the lifting of $2\beta$ is one of the  boundary components of $U_T$. By analyzing the boundaries of $U_T$ we can show that the lifting of $2\beta$ into $W$ is not homologous to zero and this fact will imply that $\beta$ is a generator of $H_1(\overline{Z};\mathbb{Z})$.

In this paper we construct four different families of numerical Campedelli surfaces with torsion $\mathbb{Z}/4\mathbb{Z}$. We will show that there is a $\mathbb{Q}$-Gorenstein smoothing of the singular surface $X$ in the main construction such that a general fiber $X_t$ of the $\mathbb{Q}$-Gorenstein smoothing of $X$ has a non-ample canonical divisor; Proposition~\ref{proposition:main-not-ample}. In \S\ref{subsection:main-modified-ample}, by modifying the configuration used in the main construction, we construct another numerical Campedelli surface with $H_1=\mathbb{Z}/4\mathbb{Z}$ whose canonical divisor is ample. Furthermore, we construct two more examples of numerical Campedelli surface with $H_1=\mathbb{Z}/4\mathbb{Z}$ using different configurations in \S\ref{subsection:more-not-ample}. Finally, we also construct a numerical Campedelli surface with $H_1=\mathbb{Z}/2\mathbb{Z} \oplus \mathbb{Z}/2\mathbb{Z}$ using a similar method in Appendix.

\subsubsection*{Note}

Few months after this paper was announced, Frapporti~\cite{Frapporti} constructed a numerical Campedelli surface with $\pi_1=\mathbb{Z}/4\mathbb{Z}$ by a completely different method.

\subsection*{Acknowledgements}

The authors would like to thank Professor Miles Reid for suggesting this problem, and Professor JongHae Keum and Professor Ki-Heon Yun for helpful discussion during the work. The authors also wish to thank Professor Roberto Pignatelli for explaining how to compute the algebraic fundamental groups of numerical Campedelli surfaces. Heesang Park was supported by Basic Science Research Program through the National Research Foundation of Korea (NRF) grant funded by the Korean Government (2011-0012111). Jongil Park was supported by the National Research Foundation of Korea (NRF) grant funded by the Korean Government (2010-0019516). He also holds a joint appointment at KIAS and in the Research Institute of Mathematics, SNU. Dongsoo Shin was supported by Basic Science Research Program through the National Research Foundation of Korea (NRF) grant funded by the Korean Government (2010-0002678).

\section{An Enriques surface}

We start with an elliptic Enriques surface constructed by Kondo~\cite[Example~V]{Kondo}. We briefly summarize the construction for the convenience of the reader. We will follow the notations of Kondo~\cite{Kondo}.

\begin{figure}[htb]
    \centering
    \includegraphics[scale=0.8]{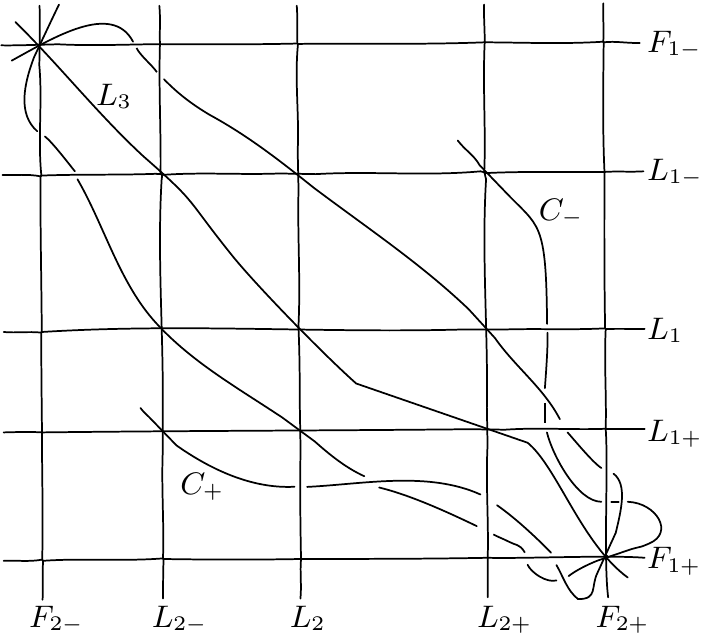}
    \caption{$\mathbb{P}^1 \times \mathbb{P}^1$; cf.\thinspace Kondo~\cite[Fig. 5.1]{Kondo}}
    \label{figure:P1timesP1}
    \end{figure}

Let $\mathbb{P}^1 \times \mathbb{P}^1=\{(u_0:u_1, v_0:v_1) \mid [u_0:u_1], [v_0:v_1] \in \mathbb{P}^1\}$. Let us consider the following smooth rational curves:
    \begin{align*}
    &C_{+}: (v_0+v_1)(u_0+u_1)=-2(u_0-u_1)v_0,\\
    &C_{-}: (v_0-v_1)(u_0-u_1)=-2(u_0+u_1)v_0,\\
    &L_{1\pm}: u_0=\pm u_1, \quad L_{2\pm}: v_0=\pm v_1,\\
    &L_1: u_0=0, \quad L_2: v_0=0, \quad L_3: u_0v_1=u_1v_0,\\
    &F_{1\pm}: u_0=\pm\frac{1}{\sqrt{-3}} u_1, \quad F_{2\pm}: v_0=\pm \frac{1}{\sqrt{-3}} v_1.
    \end{align*}
The configurations of these curves are given in Figure~\ref{figure:P1timesP1}.

\begin{figure}[htb]
    \centering
    \includegraphics[scale=0.9]{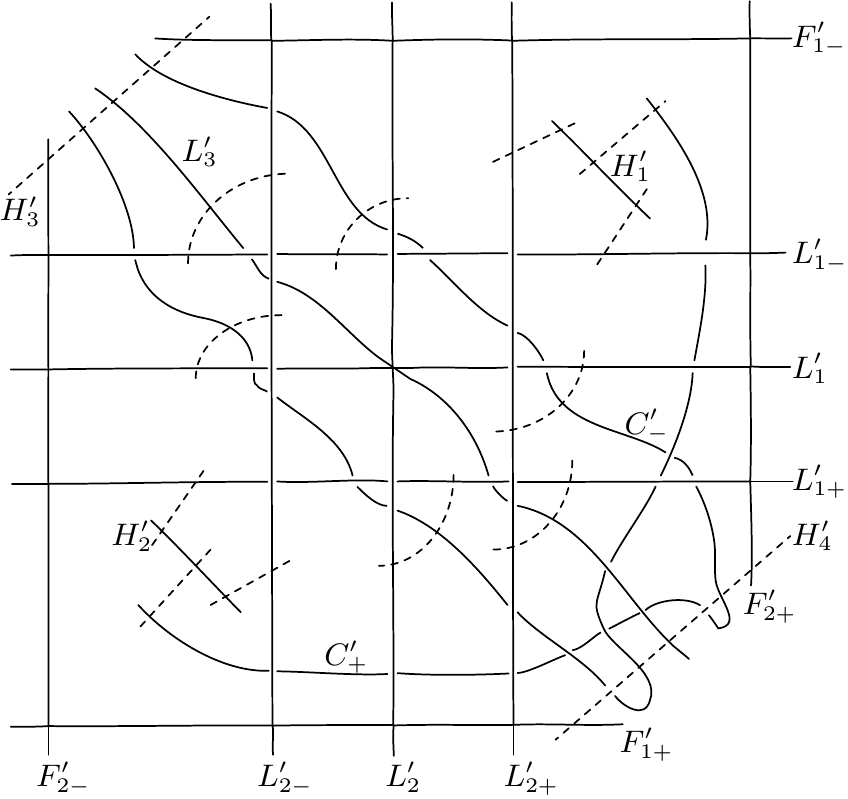}
    \caption{$R=(\mathbb{P}^1 \times \mathbb{P}^1) \sharp 16 \overline{\mathbb{P}^2}$; cf.\thinspace Kondo~\cite[Fig. 5.2]{Kondo}}
    \label{figure:P1timesP1+16P2}
    \end{figure}

Put $B=C_{+}+C_{-}+L_{1+}+L_{1-}+L_{2+}+L_{2-}$, which is a $(4,4)$-divisor in $\mathbb{P}^1 \times \mathbb{P}^1$. We first blow up $\mathbb{P}^1 \times \mathbb{P}^1$ at the $10$ singular points of $B$. Let $H_1$ and $H_2$ be the exceptional curves over $(-1:1,1:1)$ and $(1:1,-1:1)$, respectively. We blow up again at the $6$ intersection points of $H_1$, $H_2$ and the proper transforms of $L_{1\pm}$, $L_{2\pm}$, and $C_{\pm}$. Denote by $R=(\mathbb{P}^1 \times \mathbb{P}^1) \sharp 16\overline{\mathbb{P}^2}$ the blown up rational surface; Figure~\ref{figure:P1timesP1+16P2}. Let $C_{\pm}'$, $L_{1\pm}'$, $L_{2\pm}'$, $L_1'$, $L_2'$, $L_3'$, $F_{1\pm}'$, $F_{2\pm}'$, $H_1'$, $H_2'$ be the proper transforms of $C_{\pm}$, $L_{1\pm}$, $L_{2\pm}$, $L_1$, $L_2$, $L_3$, $F_{1\pm}$, $F_{2\pm}$, $H_1$, $H_2$, respectively. The curves $H_3'$ and $H_4'$ denote the proper transforms of the exceptional curves over $(-1/\sqrt{-3}:1,-1/\sqrt{-3}:1)$ and $(1/\sqrt{-3}:1,1/\sqrt{-3}:1)$, respectively. The configurations of the proper transforms of the curves and the exceptional curves in $R$ are given in Figure~\ref{figure:P1timesP1+16P2}.

\subsection{An Enriques surface}

\begin{figure}[htb]
    \centering
    \includegraphics[scale=0.9]{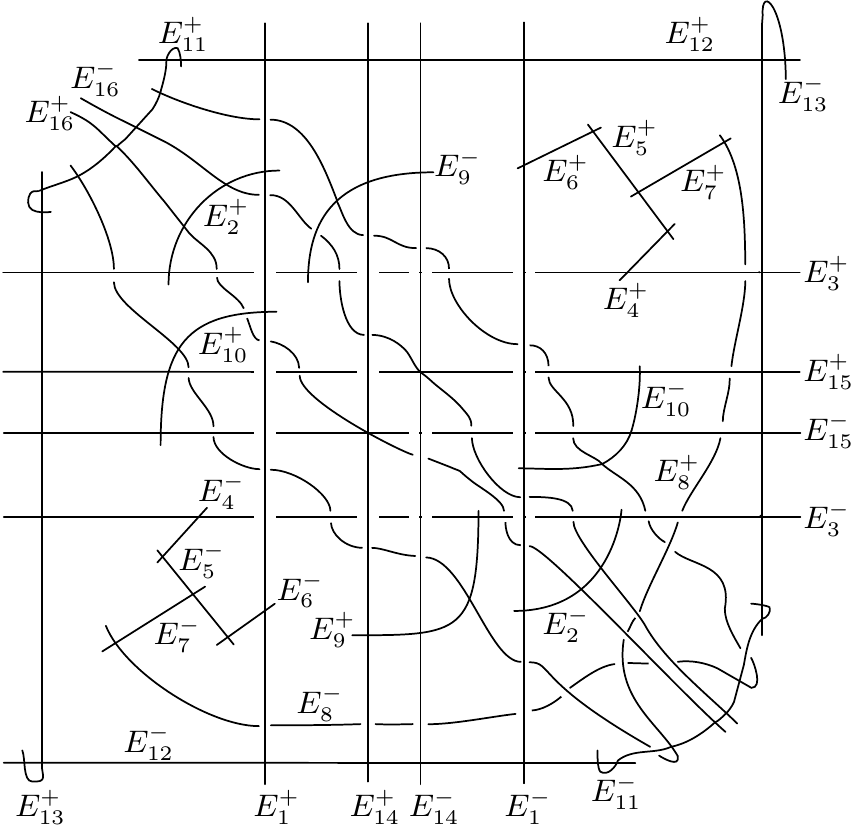}
    \caption{A K3 surface $V$; cf.\thinspace Kondo~\cite[Fig. 5.3]{Kondo}}
    \label{figure:K3}
    \end{figure}

Let $V$ be the double covering of $R$ branched along the divisor $B'=C_+'+C_-'+L_{1+}'+L_{1-}'+L_{2+}'+L_{2-}'+H_1'+H_2'$. Then $V$ is a K3 surface. Let $\phi: V \to R$ be the covering map. The configuration of the inverse images of the rational curves of $R$ by the covering $\phi$ is given in Figure~\ref{figure:K3}.

\begin{figure}[htb]
    \centering
    \includegraphics[scale=0.9]{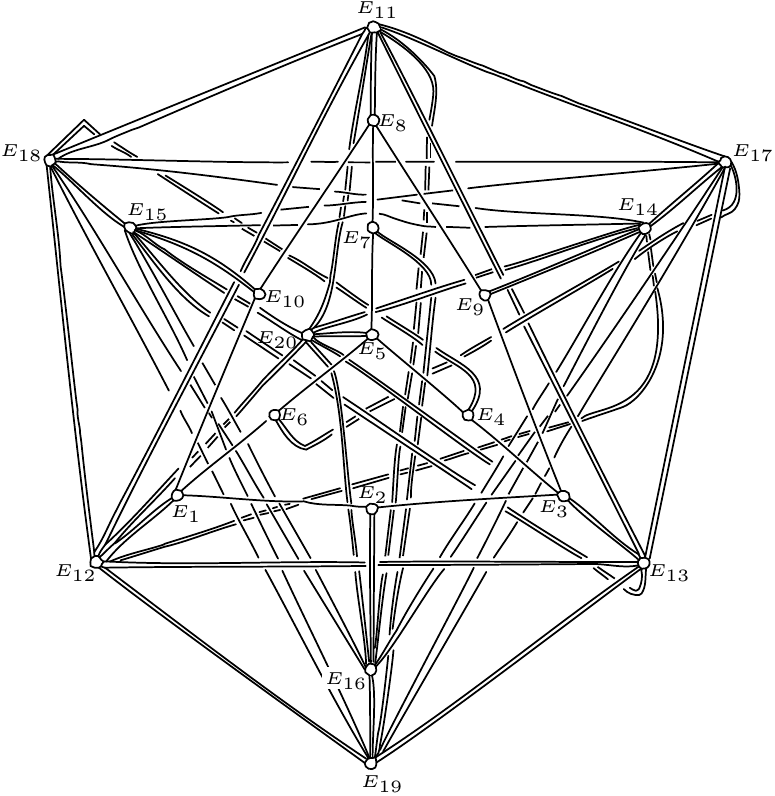}
    \caption{A dual graph of rational curves on $Y$; cf.\thinspace Kondo~\cite[Fig. 5.5]{Kondo}}
    \label{figure:dual_graph}
    \end{figure}

The K3 surface $V$ has a fixed point free involution $\sigma: V \to V$ defined by the composition of the involution induced by the map $\nu: \mathbb{P}^1 \times \mathbb{P}^1 \to \mathbb{P}^1 \times \mathbb{P}^1$, $(u_0:u_1, v_0:v_1) \mapsto (-u_0:u_1, -v_0:v_1)$ and the covering involution $\tau: V \to V$ induced by the covering $\phi: V \to R$.  Then the quotient $Y=V/\langle \sigma \rangle$ is an Enriques surface. Let $\pi : V \to Y$ be the unramified double covering. Since $\sigma$ acts on the rational curves $E_i^\pm$ as $\sigma(E_i^+)=E_i^-$, we have $\pi(E_i^+)=\pi(E_i^-)$. Set $E_i=\pi(E_i^+)=\pi(E_i^-)$. According to Kondo~\cite[Example V]{Kondo}, there are exactly 20 rational curves in $Y$. The dual graph of all rational curves on $Y$ is given in Figure~\ref{figure:dual_graph}.

\subsection{An elliptic fibration}

\begin{figure}[htb]
    \centering
    \includegraphics[scale=0.9]{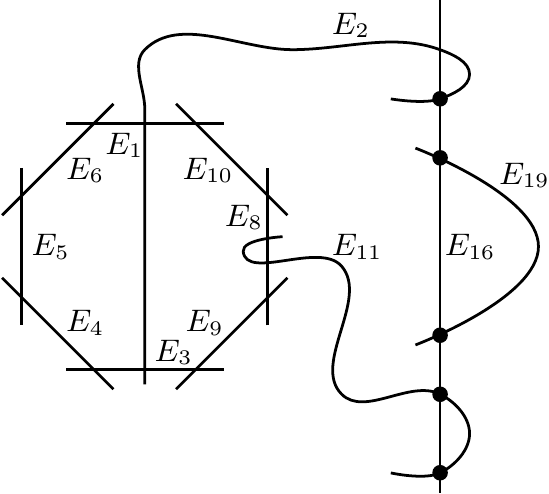}
    \caption{An Elliptic fibration on $Y$}
    \label{figure:Y}
    \end{figure}

In the Enriques surface $Y$ the linear system  $|E_1+E_6+E_5+E_4+E_3+E_9+E_8+E_{10}|=|E_{16}+E_{19}|$ defines an elliptic fibration with one $I_8$-singular fiber consisting of $E_1+E_6+E_5+E_4+E_3+E_9+E_8+E_{10}$ and one $I_2$-singular fiber consisting of $E_{16}+E_{19}$. The two rational curves $E_2$ and $E_{11}$ are bisections of the elliptic fibration; Figure~\ref{figure:Y}. Note that $E_2$ intersects transversely $E_{16}$ at two different points. Also $E_{11}$ intersects $E_8$ and $E_{16}$ at two different points, respectively; cf.\thinspace Figure~\ref{figure:P1timesP1+16P2}.

\begin{remark}\label{remark:E_19}
All the rational curves except $E_{19}$ in the singular fibers $I_8+I_2$ of the elliptic fibration of $Y$ are originated from the curves given in Figure~\ref{figure:P1timesP1+16P2}. The curves $E_3$ and $E_{16}$ were the rational curves $L_{1\pm}$ and $L_3$ (respectively) in $\mathbb{P}^1 \times \mathbb{P}^1$ and the other curves are exceptional curves contained in $R=(\mathbb{P}^1 \times \mathbb{P}^1) \sharp 16 \overline{\mathbb{P}^2}$; cf.\thinspace Figure~\ref{figure:P1timesP1}.

On the other hand the rational curve $E_{19}$ is induced from the curve $G: u_0v_0+u_1v_1=0$ in $\mathbb{P}^1 \times \mathbb{P}^1$. The curve $G$ intersects with the other curves $C_{\pm}$, $L_{1\pm}$, $L_{2\pm}$, $L_1$, $L_2$, $F_{1\pm}$, $F_{2\pm}$ as follows:
    {\allowdisplaybreaks
    \begin{align*}
    &G \cap C_{\pm} = G \cap L_{1\pm} = G \cap L_{2\mp} = \{(\pm 1:1, \mp 1:1)\}, \\
    &G \cap L_1 = \{(0:1, 1:0)\}, \ G \cap L_2 = \{(1:0, 0:1)\},\\
    &G \cap L_3 = \{(\sqrt{-1}:1,\sqrt{-1}:1), (-\sqrt{-1}:1,-\sqrt{-1}:1)\},\\
    &G \cap F_{1\pm} = \left\{\left(\pm \frac{1}{\sqrt{-3}}:1, \mp \sqrt{-3}:1\right)\right\},\\
    &G \cap F_{2\pm} = \left\{\left(\mp \sqrt{-3}:1, \pm \frac{1}{\sqrt{-3}}:1\right)\right\}.
    \end{align*}
    }
Note that $G$ is tangent to $C_+$ and $C_-$ at the intersection points, respectively. Since $G$ intersects with $L_3$ at two different points, their images $E_{16}$ and $E_{19}$ in the Enriques surface $Y$ form an $I_2$-singular fiber.
\end{remark}

\section{Main construction}\label{section:main-construction}

In this section we first construct a singular surface $X$ with three permissible singularities. We then apply a $\mathbb{Q}$-Gorenstein smoothing theory to the singular surface $X$ so that we obtain the desired surface.

\begin{figure}[htb]
    \centering
    \includegraphics[scale=0.9]{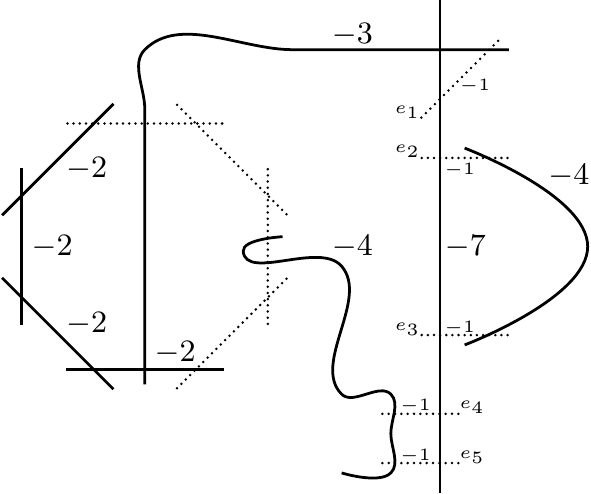}
    \caption{$Z = Y \sharp 5\overline{\mathbb{P}^2}$}
    \label{figure:Z}
    \end{figure}

We blow up at the five marked points $\bullet$ on the Enriques surface $Y$; cf.\thinspace Figure~\ref{figure:Y}. Then we get a surface $Z=Y \sharp 5\overline{\mathbb{P}^2}$; Figure~\ref{figure:Z}. We denote the five exceptional curves by $e_1, \dotsc, e_5$. There exist three disjoint linear chains of $\mathbb{P}^1$'s in $Z$, which are denoted by the following dual graphs:
    \begin{equation*}
    \begin{aligned}
    &C_1: \udc{-7}{u_1}-\udc{-3}{u_2}-\udc{-2}{u_3}-\udc{-2}{u_4}-\udc{-2}{u_5}-\udc{-2}{u_6},\\
    &C_2: \udc{-4}{u_7}, \quad C_3: \udc{-4}{u_8},
    \end{aligned}
    \end{equation*}
where $u_i$ denotes the corresponding rational curves. The configuration $C_1$ consists of the proper transforms of $E_{16}$, $E_2$, $E_3$, $E_4$, $E_5$, $E_6$, $C_2$ consists of the proper transform of $E_{11}$, and $C_3$ consists of the proper transform of $E_{19}$.

We contract these three chains of $\mathbb{P}^1$'s from the surface $Z$ so that it produces a normal projective surface $X$ with three singularities of class $T$. In Section~\ref{section:global-Q} we will show that the singular surface $X$ has a global $\mathbb{Q}$-Gorenstein smoothing; Theorem~\ref{theorem:Q-Gorenstein}. Let $X_t$ be a general fiber of a $\mathbb{Q}$-Gorenstein smoothing of the singular surface $X$.

\begin{theorem}\label{theorem:properties_of_X_t}
A general fiber $X_t$ is a minimal complex surface of general type with $p_g=0$, $K^2=2$ and $H_1(X_t;\mathbb{Z})=\mathbb{Z}/4\mathbb{Z}$.
\end{theorem}

\begin{proof}
Since $X$ is a singular surface with $p_g=0$ and $K^2=2$, by applying general results of complex surface theory and $\mathbb{Q}$-Gorenstein smoothing theory, one may conclude that a general fiber $X_t$ is a complex surface with $p_g=0$ and $K^2=2$.

We now prove the minimality of $X_t$. Let $f: Z \to X$ be the contraction map and let $h: Z \to Y$ be the blowing-up. Then we have
    \begin{align*}
    K_Z &= f^{\ast}{K_X} - \left( \frac{5}{6} u_1 + \frac{5}{6} u_2 + \frac{4}{6} u_3 + \frac{3}{6} u_4 + \frac{2}{6} u_5 + \frac{1}{6} u_6 \right) - \frac{1}{2} u_7 - \frac{1}{2} u_8,\\
    K_Z &= h^{\ast}{K_Y} + e_1 + e_2 + e_3 + e_4 + e_5,
    \end{align*}
Since $K_Y$ is numerically trivial, by combining these relations, we get
    \begin{equation*}
    \begin{split}
    f^{\ast}{K_X} \equiv &\frac{5}{6} u_1 + \frac{5}{6} u_2 + \frac{4}{6} u_3 + \frac{3}{6} u_4 + \frac{2}{6} u_5 + \frac{1}{6} u_6 + \frac{1}{2} u_7 + \frac{1}{2} u_8 \\
    &+ e_1 + e_2 + e_3 + e_4 + e_5,
    \end{split}
    \end{equation*}
where $\equiv$ denotes the numerical equivalence. Since all the coefficients are positive in the expression of $f^{\ast}{K_X}$, the $\mathbb{Q}$-divisor $f^{\ast}{K_X}$ is nef if $f^{\ast}{K_X} \cdot e_i \ge 0$ for all $i=1, \dotsc, 5$. In fact we have $f^{\ast}{K_X} \cdot e_1 = \frac{2}{3}$, $f^{\ast}{K_X} \cdot e_2 = \frac{1}{3}$, $f^{\ast}{K_X} \cdot e_3 = \frac{1}{3}$, $f^{\ast}{K_X} \cdot e_4 = \frac{1}{3}$, $f^{\ast}{K_X} \cdot e_5 = \frac{1}{3}$. Therefore $f^{\ast}{K_X}$ is nef; hence $K_X$ is also nef. Let $\pi: \chi \to \Delta$ be a $\mathbb{Q}$-Gorenstein smoothing of $X$. Since the  $\mathbb{Q}$-Cartier divisor $K_{\chi/\Delta}$ is $\pi$-big over $\Delta$ and  $\pi$-nef at the point $0$, the nefness of $K_{X_t}$ is also obtained by shrinking $\Delta$ if it is necessary; cf.\thinspace Nakayama~\cite{Nakayama}. Therefore a general fiber $X_t$ is minimal.

Since $K_{X_t}$ is nef, the Kodaira dimension of $X_t$ is nonnegative. By Enriques-Kodaira classification, all minimal surfaces with Kodaira dimension $0$ or $1$ have $K^2=0$. But, $K_{X_t}^2=2$; hence, $X_t$ is of general type.

Finally we will prove that $H_1(X_t;\mathbb{Z})=\mathbb{Z}/4\mathbb{Z}$ in Proposition~\ref{proposition:H_1=Z_4}.
\end{proof}

\begin{proposition}\label{proposition:main-not-ample}
There is a $\mathbb{Q}$-Gorenstein smoothing of $X$ such that the canonical divisor $K_{X_t}$ of a general fiber $X_t$ is \emph{not} ample.
\end{proposition}

\begin{proof}
The proper transform $\widetilde{E}_{10}$ of $E_{10}$ is a ($-2$)-curve in $Z$; cf.~Figure~\ref{figure:Y}, \ref{figure:Z}. We now construct a $\mathbb{Q}$-Gorenstein smoothing of the singular surface $X$ which can be extended to the deformation of the pair $(X, \widetilde{E}_{10})$.

We contract once more $\widetilde{E}_{10}$ from $X$ so that we obtain a singular surface $X'$ with four singular points of class $T$: Three of them are the original singular points on $X$, say, $p_1, p_2, p_3$, and the new one is a rational double point, denoted by $q$. By applying a similar method as in the proof of Theorem~\ref{theorem:Q-Gorenstein}, it is not difficult to show that the local-to-global obstruction space of the singular surface $X'$ vanishes also. Therefore every local deformations of the four singularities can be globalized. Consider the local deformation of $X'$ consisting of $\mathbb{Q}$-Gorenstein smoothings of the original singularities $p_1, p_2, p_3$ and the trivial deformation of the rational double point $q$. Then the local deformation is globalized and we obtain a family $\mathcal{X}' \to \Delta$ over a small disk $\Delta$ with the central fiber $X'_0 = X'$. Note that every fiber $\mathcal{X}'_t$ has a rational double point as its unique singularity. By shrinking $\Delta$ if necessary, we resolve simultaneously the singularity of each fiber $X'_t$. We then obtain a family $\mathcal{X}'' \to \Delta$ such that every fiber $X''_t$ has a ($-2$)-curve. Note that the central fiber of the family $\mathcal{X}'' \to \Delta$ is the singular surface $X$. Therefore we construct a $\mathbb{Q}$-Gorenstein smoothing of $X$ that is extended to a deformation of the pair $(X, \widetilde{E}_{10})$.

By applying the same method in the proof of the above Theorem~\ref{theorem:properties_of_X_t}, we can show that a general fiber $X''_t$ is a minimal complex surface of general type with $p_g=0$, $K^2=2$ and $H_1=\mathbb{Z}/4\mathbb{Z}$. But, since a general fiber $X''_t$ contains a ($-2$)-curve, the canonical divisor $K_{X''_t}$ is not ample.
\end{proof}

\begin{remark}
We don't know that the canonical divisor $K_{X_t}$ of a general fiber $X_t$ is \emph{not} ample for any $\mathbb{Q}$-Gorenstein smoothing of $X$. In fact we will construct an example with the ample canonical divisor in \S\ref{subsection:main-modified-ample} by modifying the main example in the above Theorem~\ref{theorem:properties_of_X_t}.
\end{remark}

\section{Existence of a global $\mathbb{Q}$-Gorenstein smoothing}
\label{section:global-Q}

This section is devoted to a proof of the following theorem.

\begin{theorem}\label{theorem:Q-Gorenstein}
The singular surface $X$ has a global $\mathbb{Q}$-Gorenstein smoothing.
\end{theorem}

The following proposition tells us a sufficient condition for the existence of a global $\mathbb{Q}$-Gorenstein smoothing of $X$.

\begin{proposition}[Y. Lee-J. Park~\cite{Lee-Park-K^2=2}]\label{proposition:sufficient_H^2=0}
Let $X$ be a normal projective surface with singularities of class $T$. Let $f: \widetilde{X} \to X$ be the minimal resolution and let $A$ be the reduced exceptional divisor. If $H^2(T_{\widetilde{X}}(-\log{A}))=0$, then there is a global $\mathbb{Q}$-Gorenstein smoothing of $X$.
\end{proposition}

Since the contraction map $f: Z \to X$ is the minimal resolution of the singular surface $X$, the existence of a global $\mathbb{Q}$-Gorenstein smoothing of $X$ follows from the vanishing of the cohomology $H^2(T_Z(-\log{A}))$, where $A = u_1 + \dotsb + u_8$ is the divisor on $Z$ consisting of the contracted rational curves. On the one hand the vanishing $H^2(T_Z(-\log{A}))=0$ is preserved under blowing-downs.

\begin{proposition}[Flenner-Zaidenberg~\cite{Flenner-Zaidenberg}]\label{proposition:Flenner-Zaidenberg}
Let $S$ be a nonsingular surface and let $A$ be a simple normal crossing divisor in $S$. Let $f : S' \to S$ be a blowing up of $S$ at a point p of $A$. Set $A'=f^{-1}(A)_{red}$. Then $h^2(T_{S'}(-\log{A'}))=h^2(T_S(-\log{A}))$.
\end{proposition}

Let
    \[D=E_2+E_3+E_4+E_5+E_6+E_{11}+E_{16}+E_{19} \in \divisor(Y)\]
be the divisor in the Enriques surface $Y$ consisting of the rational curves whose proper transforms in the blown-up Enriques surface $Z$ are contracted. Note that $f^{-1}(D)_{\text{red}} = A + e_1+e_2+e_3+e_4+e_5$. Therefore if $h^2(T_Y(-\log{D}))=0$ then $h^2(T_Z(-\log{A}))=0$ by Proposition~\ref{proposition:Flenner-Zaidenberg}. Thus Theorem~\ref{theorem:Q-Gorenstein} follows from the following proposition.

\begin{proposition}\label{proposition:H^2=0}
$H^2(T_Y(-\log{D}))=H^0(\Omega_Y(\log{D})(K_Y))=0$.
\end{proposition}

In order to prove Proposition~\ref{proposition:H^2=0}, we follow a similar strategy as in J. Keum-Y. Lee-H. Park~\cite{Keum-Lee-Park}. That is, we consider a lifting of the divisor $D$ on the Enriques surface $Y$ into the K3 surface $V$ and then we use the push-forward map of the double covering $V \to Y$ for proving that $H^0(\Omega_Y(\log{D})(K_Y))=0$.

\begin{proof}[Proof of Proposition~\ref{proposition:H^2=0}]
Let $\pi: V \to Y$ be the unramified double covering from the K3 surface $V$ to the Enriques surface $Y$. The K3 surface $V$ has two $I_8$-singular fibers and two $I_2$-singular fibers and four sections induced from the elliptic fibration structure of $Y$; Figure~\ref{figure:V}. Note that $\pi(E_i^+)=\pi(E_i^-)=E_i$. We choose a divisor $\Delta \in \divisor(V)$ of the form
    \[\Delta=E_2^++E_2^-+E_3^-+E_4^-+E_5^-+E_6^-+E_{11}^++E_{11}^-+E_{16}^++E_{19}^+\]
so that $\Delta \leq \pi^*{D}$ and $\pi_*{\Delta}=D$; Figure~\ref{figure:V-H2}.

\begin{figure}[htb]
\centering
\includegraphics[scale=0.9]{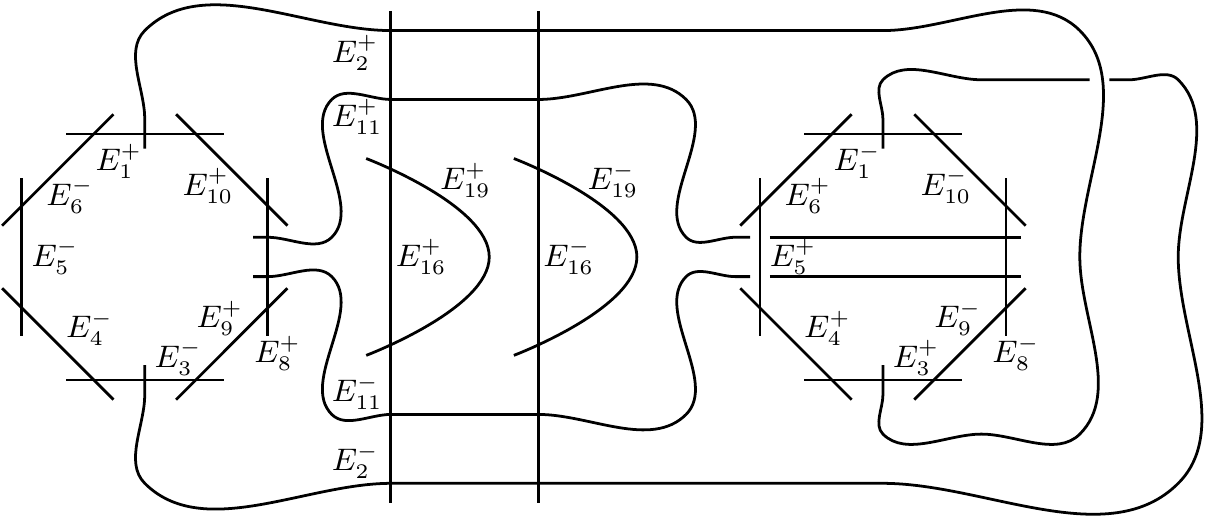}
\caption{A K3 surface $V$}
\label{figure:V}
\end{figure}

\begin{figure}[htb]
\centering
\includegraphics[scale=0.9]{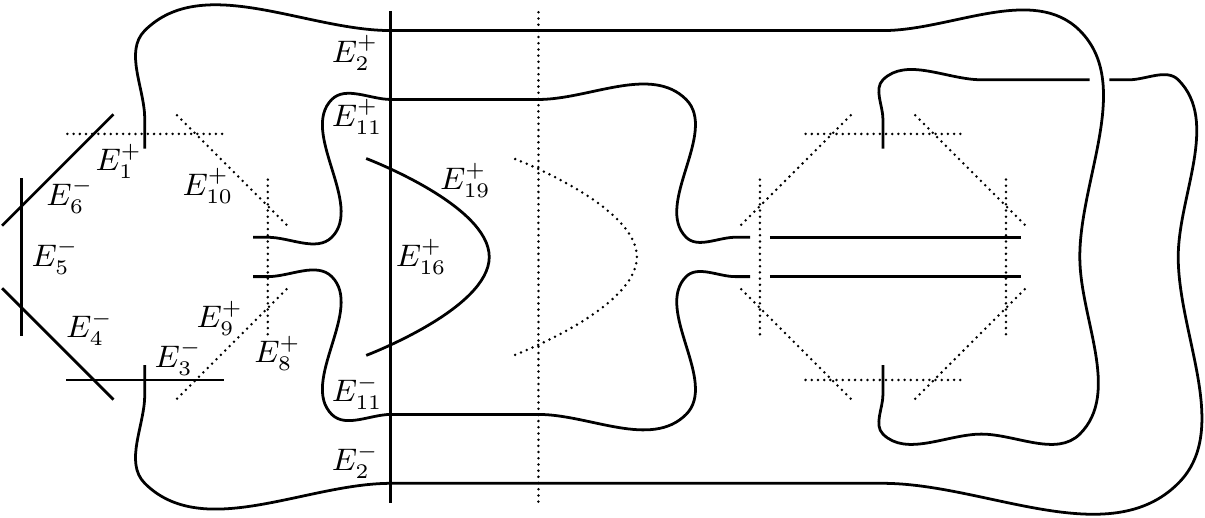}
\caption{The divisor $\Delta$}
\label{figure:V-H2}
\end{figure}

We have an exact sequence
    \begin{equation*}
    0 \to \Omega_V \to \Omega_V(\log{\Delta}) \to \bigoplus_i \sheaf{O}_{\Delta_i} \to 0,
    \end{equation*}
where $\Delta_i$'s are irreducible components of the divisor $\Delta$. In the long exact sequence
    \begin{equation*}
    0 \to H^0(\Omega_V) \to H^0(\Omega_V(\log{\Delta})) \to \bigoplus_i H^0(\sheaf{O}_{\Delta_i}) \xrightarrow{\delta} H^1(\Omega_V^1)\to \cdots
    \end{equation*}
the connecting homomorphism $\delta: \bigoplus_i H^0( \sheaf{O}_{\Delta_i}) \to H^1(\Omega_V^1)$ is the first Chern class map. Since the intersection matrix of the irreducible components of $\Delta$ is invertible, their images by the first Chern class map $\delta$ are linearly independent. Hence $\delta$ is injective. Furthermore, since $V$ is a K3 surface, $H^0(\Omega_V)=0$. Therefore we have
    \[H^0(V, \Omega_V(\log \Delta))=0.\]
By the choice of $\Delta$, we have $\Omega_Y(\log D)\subset \pi_*(\Omega_V(\log \Delta))$. On the other hand, since $0=K_V=\pi^*K_Y$, it follows by the projection formula that
    \[\pi_*(\Omega_V(\log \Delta)) = \pi_*(\Omega_V(\log \Delta)(K_V))=\pi_*(\Omega_V(\log \Delta))\otimes K_Y.\]
Therefore we have
    \begin{equation*}
    \begin{split}
    H^0(Y,\Omega_Y(\log D)(K_Y)) &\subset H^0(Y, \pi_*(\Omega_V(\log \Delta))(K_Y)) \\
    &=H^0(Y,\pi_*(\Omega_V(\log \Delta))) \\
    &= H^0(V, \Omega_V(\log \Delta)) \\
    &=0. \qedhere
    \end{split}
    \end{equation*}
\end{proof}

\section{Proof of $H_1(X_t; \mathbb{Z})=\mathbb{Z}/4\mathbb{Z}$}
\label{section:Z_4}

In this section we calculate the first homology group of a general fiber $X_t$
of a $\mathbb{Q}$-Gorenstein smoothing of $X$, which is the key part of this article.

\begin{theorem}\label{proposition:H_1=Z_4}
Let $X_t$ be a general fiber of a $\mathbb{Q}$-Gorenstein smoothing of the singular surface $X$. Then $H_1(X_t; \mathbb{Z})=\mathbb{Z}/4\mathbb{Z}$.
\end{theorem}

We denote by
    \[C_1: \udc{-7}{u_1}-\udc{-3}{u_2}-\udc{-2}{u_3}-\udc{-2}{u_4}-\udc{-2}{u_5}-\udc{-2}{u_6}, \quad C_2: \udc{-4}{u_7}, \quad C_3: \udc{-4}{u_8}
    \]
the plumbings of the linear chains of the rational curves. Let $\overline{Z}$ be a smooth  $4$-manifold obtained from $Z$ by replacing three plumbings $C_1$, $C_2$, $C_3$ with the corresponding Milnor fibers $M_1$, $M_2$, $M_3$, respectively. Since a general fiber $X_t$ is diffeomorphic to the rational blow-down $4$-manifold $\overline{Z}$ by Milnor fiber theory, we have $H_1(X_t; \mathbb{Z})=H_1(\overline{Z}; \mathbb{Z})$. Hence it suffices to show that
    \[H_1(\overline{Z}; \mathbb{Z}) = \mathbb{Z}/4\mathbb{Z}.\]
We will prove it in Proposition~\ref{proposition:H_1(Zbar)=Z_4}. We start with decomposing the surface $Z$ into
    \[Z=Z_0 \cup (C_1 \cup C_2 \cup C_3).\]
Then the $4$-manifold $\overline{Z}$ can be decomposed into
    \[\overline{Z}=Z_0 \cup (M_1 \cup M_2 \cup M_3).\]
In order to prove that $H_1(\overline{Z}; \mathbb{Z})=\mathbb{Z}/4\mathbb{Z}$, we first calculate $H_1(Z_0;\mathbb{Z})$.

\begin{proposition}\label{proposition:H_1(Z0)=Z_4}
$H_1(Z_0;\mathbb{Z})=\mathbb{Z}/4\mathbb{Z}$.
\end{proposition}

We divide the proof of Proposition~\ref{proposition:H_1(Z0)=Z_4} into the following several lemmas. Proposition~\ref{proposition:H_1(Z0)=Z_4} will be proved in Lemma~\ref{lemma:beta}.

Let $\alpha_1$, $\alpha_2$, $\alpha_3$ be the normal circles of the disk bundles $C_1$, $C_2$, $C_3$ over $u_1$, $u_7$, $u_8$, respectively; cf.\thinspace Figure~\ref{figure:Z}. We denote again by $\alpha_i$ the homology classes of $\alpha_i$ in $H_1(Z_0;\mathbb{Z})$ for convenience.

\begin{lemma}\label{lamme:2alpha=0}
$\alpha_1=\alpha_2=\alpha_3$ and $2\alpha_1=2\alpha_2=2\alpha_3=0$ in $H_1(Z_0;\mathbb{Z})$.
\end{lemma}

\begin{proof}
Refer Figure~\ref{figure:Z}. Since the exceptional curve $e_2$ intersects transversely the curves $u_1$ and $u_8$, we have
    \[\alpha_1+\alpha_3=0\]
in $H_1(Z_0;\mathbb{Z})$. Since the exceptional curve $e_4$ intersects transversely the curves $u_1$ and $u_7$, we have
    \[\alpha_1+\alpha_2=0\]
in $H_1(Z_0;\mathbb{Z})$. Since the curve $u_7$ intersects transversely the proper transform of $E_8$ at two different points, we have
    \[2\alpha_2=0.\]
Combining these relations, we have $\alpha_1=\alpha_2=\alpha_3$ and $2\alpha_1=2\alpha_2=2\alpha_3=0$ in $H_1(Z_0;\mathbb{Z})$.
\end{proof}

From now on we set $\alpha=\alpha_1=\alpha_2=\alpha_3 \in H_1(Z_0;\mathbb{Z})$.

\begin{lemma}\label{lemma:|H_1(Z_0)|=4}
$H_1(Z_0;\mathbb{Z})= \mathbb{Z}/4 \mathbb{Z}$ or $H_1(Z_0;\mathbb{Z})= \mathbb{Z}/2 \mathbb{Z} \oplus \mathbb{Z}/2 \mathbb{Z}$.
\end{lemma}

\begin{proof}
We first consider the Mayer-Vietoris sequence of a pair $(Z_0, C_{1} \cup C_{2} \cup  C_{3})$:
    \begin{multline*}\label{eqation:MV(Z,C)}
    H_1(\partial Z_0;\mathbb{Z}) \xrightarrow{k_{\ast}} H_1(Z_0;\mathbb{Z}) \oplus H_1(C_{1} \cup C_{2} \cup C_{3};\mathbb{Z}) \to H_1(Z;\mathbb{Z}) \\
    \to \widetilde{H_0}(\partial Z_0;\mathbb{Z}) \xrightarrow{\xi}
    \widetilde{H_0}(Z_0;\mathbb{Z}) \oplus \widetilde{H_0}(C_{1} \cup C_{2} \cup C_{3};\mathbb{Z}) \to 0.
    \end{multline*}
Since $\xi: \widetilde{H_0}(\partial Z_0;\mathbb{Z}) \to \widetilde{H_0}( C_{1} \cup C_{2} \cup  C_{3};\mathbb{Z})$ is isomorphic, we get the following short exact sequence
    \begin{equation}\label{equation:->H_1(Z_0)->}
    0 \to \im (k_{\ast}) \hookrightarrow H_1(Z_0;\mathbb{Z}) \to H_1(Z;\mathbb{Z}) \to 0.
    \end{equation}

We \textit{claim} that $\im (k_{\ast}) = \langle \alpha \rangle \cong \mathbb{Z}/2 \mathbb{Z}$. Then, since $\pi_1(Z)=\mathbb{Z}/2\mathbb{Z} = H_1(Z;\mathbb{Z})$, it follows that $H_1(Z_0;\mathbb{Z})= \mathbb{Z}/4 \mathbb{Z}$ or $H_1(Z_0;\mathbb{Z})= \mathbb{Z}/2 \mathbb{Z} \oplus \mathbb{Z}/2 \mathbb{Z}$.

\textit{Proof of claim}: By Lemma~\ref{lamme:2alpha=0} above, it suffices to show $\im (k_{\ast}) \neq 0$ in $H_1(Z_0;\mathbb{Z})$. Assume on the contrary that $\im (k_{\ast})=0$ and we consider an exact sequence of a pair $(Z_0, \partial{Z_0} = \partial{C_{1}} \cup \partial{C_{2}} \cup \partial{C_{3}})$:
    \begin{equation*}
    \cdots \rightarrow H_2(Z_0,\partial Z_0;\mathbb{Z}) \xrightarrow{\partial_{\ast}} H_1(\partial Z_0;\mathbb{Z}) \xrightarrow{k_{\ast}} H_1(Z_0;\mathbb{Z}) \rightarrow \cdots
    \end{equation*}
Then we conclude that the map $\partial_{\ast} : H_2(Z_0,\partial Z_0;\mathbb{Z}) \rightarrow H_1(\partial Z_0;\mathbb{Z})$ is surjective.

On the other hand, since a homomorphism
    \[j_{\ast} : H_2(Z;\mathbb{Z}) \to H_2(Z, C_1 \cup C_2 \cup C_3;\mathbb{Z})\]
induced by an inclusion is surjective and
    \[H_2(Z,C_1 \cup C_2 \cup  C_3;\mathbb{Z}) \cong H_2(Z_0,\partial C_1 \cup \partial C_2 \cup \partial C_3;\mathbb{Z})\]
by the excision principle, the space $\im({\partial}_{\ast}) = \im({\partial}_{\ast} \circ j_{\ast})$ is completely determined by the images of elements in $H_2(Z;\mathbb{Z})$. Furthermore, by choosing a suitable set $\mathcal{B}$ of homology classes expanding $H_2(Z;\mathbb{Z})$ and by computing the images of all elements in $\mathcal{B}$ under $\partial_{\ast} \circ j_{\ast}$, one can compute the space $\im({\partial}_{\ast})$.

For example, we can choose such a set
    \[\mathcal{B} =\{u_1, \dotsc, u_8, e_1, \dotsc, e_5, \widetilde{E}_1, \dotsc, \widetilde{E}_{20}\} \cup \{ \widetilde{F}_{\alpha} \mid \alpha \in \Lambda\},\]
where $u_i$'s are rational curves in $C_1\cup C_2 \cup C_3$, $e_i$'s denote all exceptional curves, $\widetilde{E}_i$'s denote the proper transforms on the blown-up space $Z = Y \sharp 5\overline{\mathbb{P}^2}$, and $\{\widetilde{F}_{\alpha} \mid \alpha \in \Lambda \}$ denotes the set of all half pencils of all possible elliptic fibration structures of $Y$. One can show that the images of $\widetilde{F}_{\alpha}$ under $\partial_{\ast} \circ j_{\ast}$ are generated by the images of $\{u_1, \dotsc, u_8, e_1, \dotsc, e_5, \widetilde{E}_1, \dotsc, \widetilde{E}_{20}\}$ under ${\partial}_{\ast} \circ j_{\ast}$. Hence, by the same argument used in the proof of Lemma 2.4 in Y. Lee-J. Park~\cite{Lee-Park-H1Z2}, we conclude that the element $(0,0,1)$ in the space
    \[H_1(\partial Z_0;\mathbb{Z}) = H_1(\partial C_1 \cup \partial C_2 \cup \partial C_3;\mathbb{Z}) \cong \mathbb{Z}/72\mathbb{Z} \oplus \mathbb{Z}/4\mathbb{Z} \oplus \mathbb{Z}/4\mathbb{Z}\]
is not contained in the space $\im({\partial}_{\ast} \circ j_{\ast}) = \im({\partial}_{\ast})$, which contradicts the surjectivity of $\partial_{\ast}$.
\end{proof}

\subsection{A special loop $\beta$ in $Z_0$}\label{subsubsection:beta}

We construct a loop $\beta$ in $H_1(Z_0;\mathbb{Z})$ which lies on a certain curve of genus $2$ in $Z$, which will be shown to be a generator of $H_1(Z_0;\mathbb{Z})$; Lemma~\ref{lemma:beta}.

    \begin{figure}[htb]
    \centering
    \subfloat[$v_0/v_1$]{\includegraphics[scale=0.9]{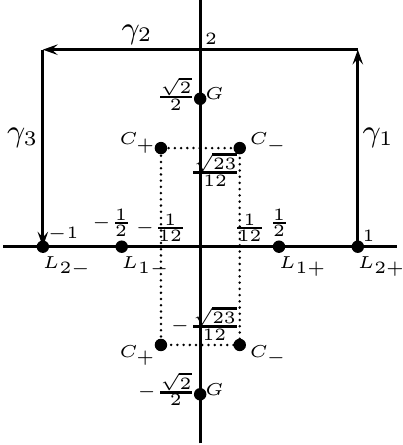}}
    \qquad
    \subfloat[The involution on $S$]{\includegraphics[scale=0.9]{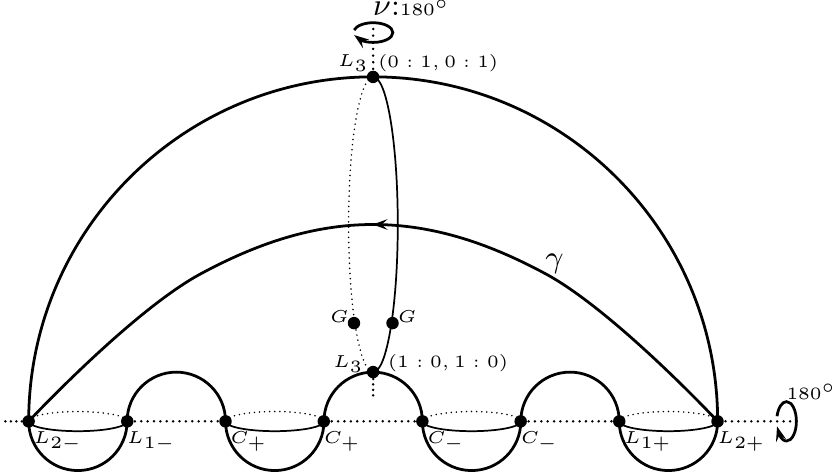}}
    \caption{A rational curve $S$ in $\mathbb{P}^1 \times \mathbb{P}^1$}
    \label{figure:S}
    \end{figure}

Set
    \[S = \{u_0v_1=2v_0u_1\} \subset \mathbb{P}^1 \times \mathbb{P}^1 = \{[u_0:u_1, v_0:v_1]\}.\]
The rational curve $S$ intersects transversally the $(4,4)$-divisor $B= C_{\pm}+L_{1\pm}+L_{2\pm}$ at eight points. In Figure~\ref{figure:S}(A) we indicate by dots $\bullet$ the $v_0/v_1$-coordinates of the intersection points of $S$ with the curves $C_{\pm}$, $L_{1\pm}$, $L_{2\pm}$ and $G$, where the horizontal axis is the real part of $v_0/v_1$ and the vertical axis is the imaginary part of $v_0/v_1$.

\begin{figure}[htb]
    \centering
    \includegraphics[scale=0.9]{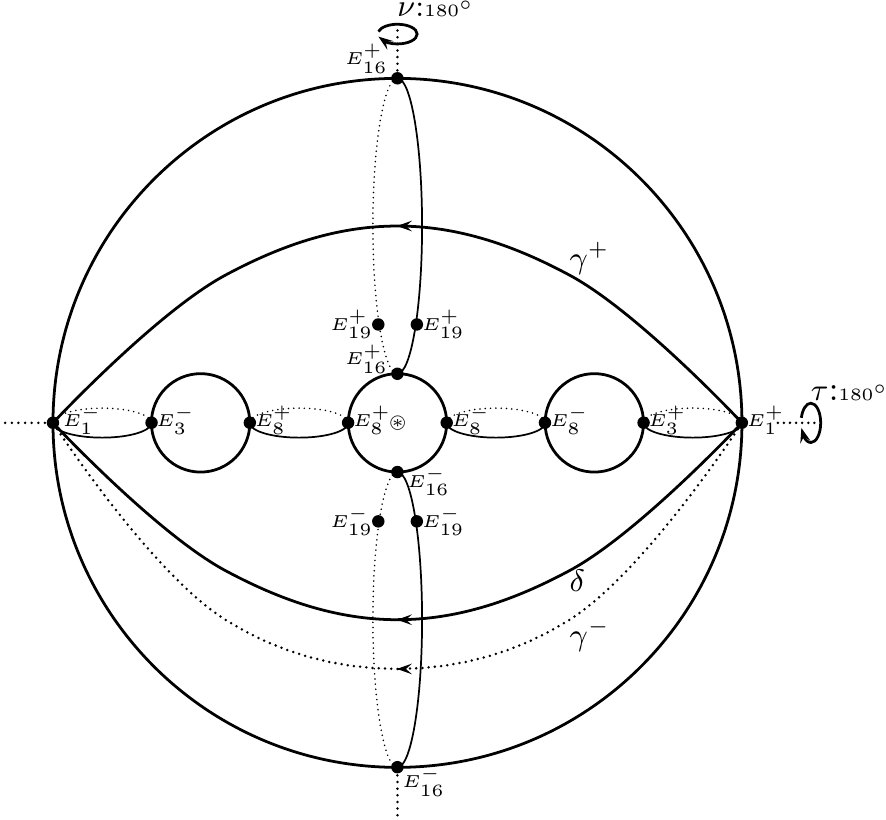}
    \caption{A curve $T$ of genus $3$ in K3 surface $V$}
    \label{figure:T}
    \end{figure}

Since $S$ does not pass through any singular points of the $(4,4)$-divisor $B$, the proper transform of $S$ in $R=(\mathbb{P}^1 \times \mathbb{P}^1) \sharp 16\overline{\mathbb{P}^2}$ intersects the branch divisor $B'=C_{\pm}'+L_{1\pm}'+L_{2\pm}'+H_1'+H_2'$ again in eight points. We denote the proper transform of $S$ again by $S$. Therefore the inverse image $T$ of the rational curve $S$ by the branched double covering $\phi: V \to R$ is a curve of genus $3$ in K3 surface $V$; Figure~\ref{figure:T}. The dots $\bullet$ indicate the intersection points with the rational curves $E_{i\pm}$.

\begin{figure}[htb]
    \centering
    \includegraphics[scale=0.9]{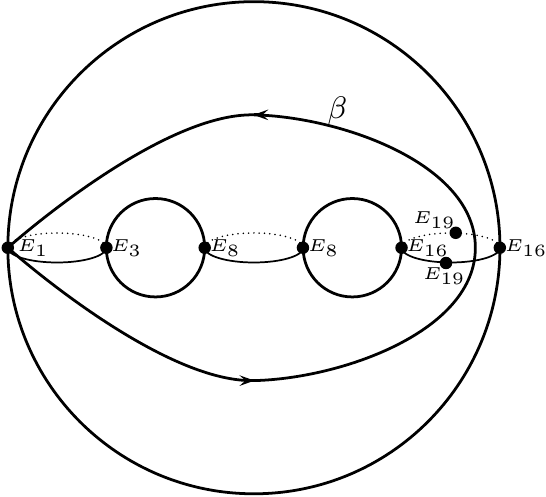}
    \caption{A curve $T/\sigma$ of genus $2$ in Enriques surface $Y$}
    \label{figure:Tmodsigma}
    \end{figure}

We denote the restriction of the double covering $\phi: V \to R$ on $T$ again by $\phi$. Hence we have a branched double covering $\phi: T \to S$. Recall that the involution $\sigma: V \to V$ that induces the Enriques surface $Y=V/\sigma$ is the composition of the extension of the involution $\nu: \mathbb{P}^1 \times \mathbb{P}^1 \to \mathbb{P}^1 \times \mathbb{P}^1$ defined by $(u_0:u_1, v_0:v_1) \mapsto (-u_0:u_1, -v_0:v_1)$ and the covering involution $\tau: V \to V$. Since the rational curve $S$ is invariant under the involution $\nu$ on $\mathbb{P}^1 \times \mathbb{P}^1$, the involution $\sigma: V \to V$ induces an involution on $T$ without fixed points. We denote the involution of $T$ again by $\sigma: T \to T$. We visualize the involutions on $T$ in Figure~\ref{figure:T}. Topologically the involution $\nu$ (or $\tau$) is the rotation around the vertical (resp. horizontal) line through the middle point  marked by $\circledast$ by an angle $180^{\circ}$. Therefore the involution $\sigma: T \to T$ is the rotation around the axis through the paper centered at the marked point $\circledast$ by an angle $180^{\circ}$. Hence the quotient $T/\sigma$ is a smooth curve of genus $2$; Figure~\ref{figure:Tmodsigma}. Let $\pi: T \to T/\sigma$ denote the quotient map.

Note that $S \cap L_{2+}'$ and $S \cap L_{2-}'$ are two of the branch points of the covering $\phi: T \to S$. Let $\gamma$ be a path lying on $S$ which connects $S \cap L_{2+}'$ and $S \cap L_{2-}'$ and does not pass through any intersection points with $C_{\pm}$, $L_{1\pm}$, $L_{2\pm}$, $L_1$, $L_2$, $L_3$, $F_{1\pm}$, $F_{2\pm}$; cf.\thinspace Figure~\ref{figure:S}(A). More precisely we define the paths $\gamma_i \subset \mathbb{P}^1 \times \mathbb{P}^1 = \{[u_0:u_1, v_0:v_1]\}$ on $S=\{u_0v_1=2v_0u_1\}$ by
    \begin{align*}
    \gamma_1(s) &= (2+2s\sqrt{-1}:1, 1+s\sqrt{-1}:1), \  0 \le s \le 2, \\
    \gamma_2(s) &= (-2s+4\sqrt{-1}:1, -s+2\sqrt{-1}:1), \ -1 \le s \le 1, \\
    \gamma_3(s) &= (-2+2(2-s)\sqrt{-1}:1, -1+(2-s)\sqrt{-1}: 1), \ 0 \le s \le 2
    \end{align*}
(cf.\thinspace Figure~\ref{figure:S}(A)), and the path $\gamma$ is defined by their union, that is,
    \begin{equation*}
    \gamma = \gamma_1 \cup \gamma_2 \cup \gamma_3.
    \end{equation*}
Note that the inverse image $\phi^{-1}(\gamma)$ on the curve $T$ in K3 surface $V$ consists of two paths $\gamma^+ \cup \gamma^-$ as in Figure~\ref{figure:T}. We define a loop $\beta$ in the curve $T/\sigma$ by the image of $\gamma^+$ under the quotient map $\pi: T \to T/\sigma$, that is,
    \begin{equation*}
    \beta=\pi(\gamma^+) \subset T/\sigma;
    \end{equation*}
cf.\thinspace Figure~\ref{figure:Tmodsigma}. Since the loop $\beta$ does not pass through any blowing-up points of the blowing-up $h : Z = Y \sharp 5\overline{\mathbb{P}^2} \to Y$, we may regard $\beta$ as a loop in $Z$.

\begin{lemma}
The loop $\beta$ is a generator of $H_1(Z;\mathbb{Z}) = \mathbb{Z}/2\mathbb{Z}$.
\end{lemma}

\begin{proof}
The path $\gamma^+$ is a lifting of $\beta \subset Y$ to K3 surface $V$, but $\gamma^+$ is not a loop. Therefore $\beta \neq 0$ in $\pi_1(Y) = \mathbb{Z}/2\mathbb{Z}$. Since $\pi_1(Y) \cong \pi_1(Z)$, we have $\beta \neq 0$ in $H_1(Z;\mathbb{Z})=\mathbb{Z}/2\mathbb{Z}$; hence it generates $H_1(Z;\mathbb{Z})$.
\end{proof}

Note that the loop $\beta$ does not pass through any intersection points of $T/\sigma$ with the rational curves consisting the configurations $C_1$, $C_2$ and $C_3$ by a choice of the path $\gamma$. Since $Z_0$ is obtained by removing $C_1$, $C_2$ and $C_3$ from $Z$, it follows that the loop $\beta$ may be considered as a loop in $Z_0$, that is, $\beta \in H_1(Z_0, \mathbb{Z})$.

\begin{lemma}\label{lemma:<alpha,beta>}
$H_1(Z_0;\mathbb{Z}) = \langle \alpha, \beta \rangle$.
\end{lemma}

\begin{proof}
In the proof of Lemma~\ref{lemma:|H_1(Z_0)|=4}, we show that $\im (k_{\ast})$ is generated by $\alpha$, and in the above lemma we show that $H_1(Z;\mathbb{Z}) = \mathbb{Z}/2\mathbb{Z}$ is generated by $\beta$. Therefore the result follows from the exact sequence \eqref{equation:->H_1(Z_0)->}.
\end{proof}

\subsection{The loop $\beta$ generates $H_1(Z_0;\mathbb{Z})$}

We will show that $H_1(Z_0;\mathbb{Z})$ is generated only by $\beta$ in Lemma~\ref{lemma:beta}, which completes the proof of $H_1(Z_0;\mathbb{Z}) = \mathbb{Z}/4\mathbb{Z}$.

\begin{figure}[htb]
    \centering
    \includegraphics[scale=0.9]{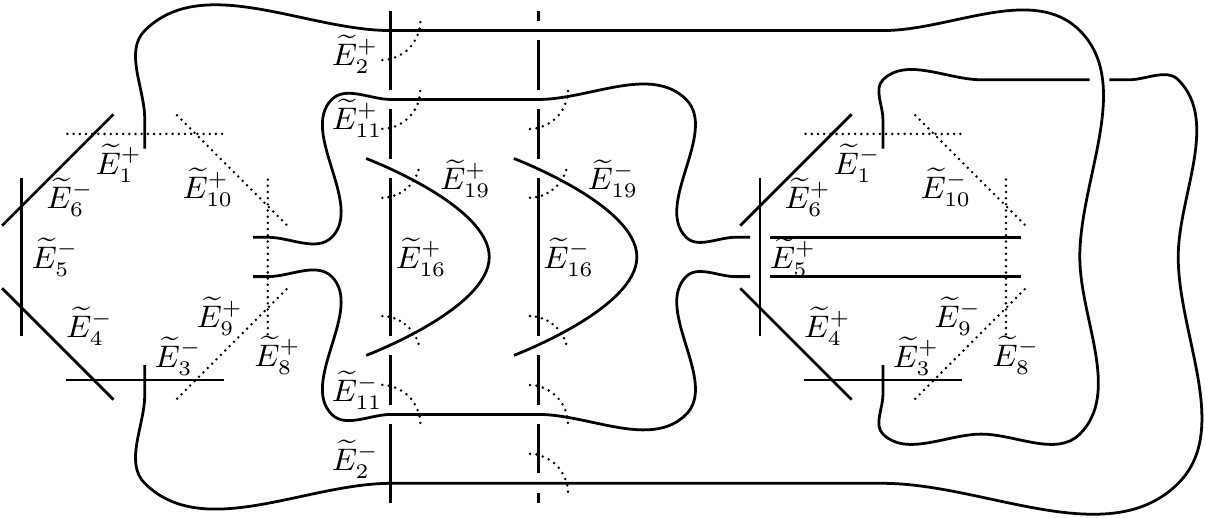}
    \caption{$W = V\sharp 10\overline{\mathbb{P}^2}$}
    \label{figure:W}
    \end{figure}

We blow up the K3 surface $V$ at the corresponding two points when we blow up the Enriques surface $Y$ at each marked point in order to construct the surface $Z$; cf.\thinspace Figure~\ref{figure:Z}. Then we get an unramified double covering $W=V\sharp 10\overline{\mathbb{P}^2} \to Z=Y\sharp 5\overline{\mathbb{P}^2}$ that induces the double covering $\pi: T \to T/\sigma$, where $T$ and $T/\sigma$ are considered as the proper transforms of $T$ and $T/\sigma$, respectively. We denote the covering $W \to Z$ again by $\pi: W \to Z$. Let us denote by $\widetilde{E}_{i}^{\pm}$ the proper transform of $E_{i}^{\pm}$. Note that the covering $W$ has six configurations of linear chains of $\mathbb{P}^1$'s:
    \begin{align*}
    &C_1^+: \udc{-7}{\widetilde{E}_{16}^+}-\udc{-3}{\widetilde{E}_2^-}-\udc{-2}{\widetilde{E}_3^-}-
    \udc{-2}{\widetilde{E}_4^-}-\udc{-2}{\widetilde{E}_5^-}-\udc{-2}{\widetilde{E}_6^-},\\
    &C_1^-: \udc{-7}{\widetilde{E}_{16}^-}-\udc{-3}{\widetilde{E}_2^+}-\udc{-2}{\widetilde{E}_3^+}-
    \udc{-2}{\widetilde{E}_4^+}-\udc{-2}{\widetilde{E}_5^+}-\udc{-2}{\widetilde{E}_6^+},\\
    &C_2^+: \udc{-4}{\widetilde{E}_{11}^+}, \quad C_2^-: \udc{-4}{\widetilde{E}_{11}^-},\\
    &C_3^+: \udc{-4}{\widetilde{E}_{19}^+}, \quad C_3^-: \udc{-4}{\widetilde{E}_{19}^-},
    \end{align*}
which are the inverse images of the configurations $C_1$, $C_2$, $C_3$ in $Z$, respectively; Figure~\ref{figure:W}. We decompose
    \[W=W_0 \cup \{C_1^+ \cup C_1^- \cup C_2^+ \cup C_2^- \cup C_3^+ \cup C_3^- \}.\]
Then the restriction $W_0 \to Z_0$ of the covering $\pi: W \to Z$ is also an unramified double covering. We denote the covering $W_0 \to Z_0$ also by $\pi: W_0 \to Z_0$ for convenience.

\begin{lemma}\label{lemma:2beta=gamma^*+alpha}
Let $\widetilde{2\beta}$ be the lifting of $2\beta \subset Z_0$ to $W_0$ by the covering map $\pi: W_0 \to Z_0$. Set $\gamma^{\ast} = \gamma^+-\gamma^-$ and let $\overline{\alpha}$ be a normal circle of $\widetilde{E}_{16}^-$ lying in $W_0$. Then $\widetilde{2\beta}$ is homologous to $\gamma^{\ast}+\overline{\alpha}$ in $W_0$.
\end{lemma}

\begin{figure}[htb]
    \centering
    \includegraphics[scale=0.9]{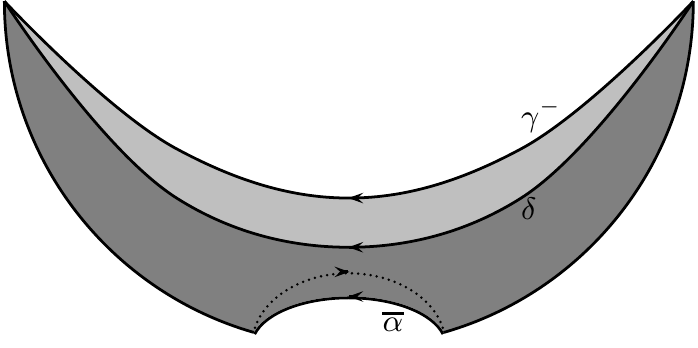}
    \caption{A lifting $\widetilde{2\beta}$}
    \label{figure:2beta}
    \end{figure}

\begin{proof}
In Figure~\ref{figure:T}, it follows that $\widetilde{2\beta}$ is homologous to $\gamma^+ - \delta$. On the other hand the curve $E_{16}^-$ intersects $T$ at two points and one of them lies in the south pole; cf.\thinspace Figure~\ref{figure:T}. The proper transform $\widetilde{E}_{16}^-$ is contained in $C_1^-$ and the configuration $C_1^-$ will be removed from $W$ to obtain in $W_0$. Since the homology class of a normal circle of the disk bundle over $\widetilde{E}_{16}^-$ is $\overline{\alpha}$, we have the following picture in $W_0$: Figure~\ref{figure:2beta}. Hence $\delta - \gamma^- + \overline{\alpha}=0$ in $H_1(W_0;\mathbb{Z})$. Therefore $\widetilde{2\beta}$ is homologous to $\gamma^+-\gamma^- + \overline{\alpha}$ in $W_0$.
\end{proof}

Since the blown-up K3 surface $W=V\sharp 10\overline{\mathbb{P}^2}$ is simply connected, we can apply the same argument used in the proof of Lemma 2.4 in Y. Lee-J. Park~\cite{Lee-Park-H1Z2} so that we can prove:

\begin{lemma}
$H_1(W_0;\mathbb{Z})\cong \mathbb{Z}/2\mathbb{Z}\cong \langle \overline{\alpha} \rangle$.
\end{lemma}

The loop $\gamma^{\ast}$ in Lemma~\ref{lemma:2beta=gamma^*+alpha} vanishes in $H_1(W_0;\mathbb{Z})$:

\begin{lemma}\label{lemma:gamma^*=0}
$\gamma^{\ast}=0$ in $H_1(W_0;\mathbb{Z})$.
\end{lemma}

\begin{figure}[htb]
    \centering
    \subfloat[$v_0/v_1$]{\includegraphics[scale=0.9]{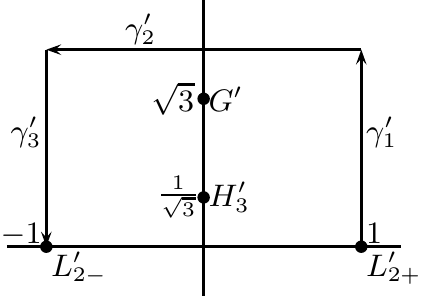}}
    \qquad
    \subfloat[the involution on $F_{1-}'$]{\includegraphics[scale=0.9]{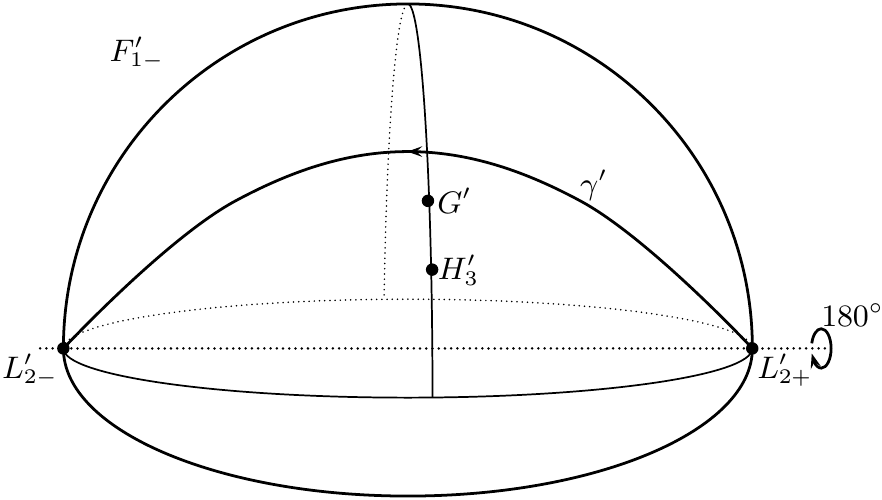}}
    \caption{$F_{1-}'$}
    \label{figure:F1-}
    \end{figure}

\begin{figure}[htb]
    \centering
    \subfloat[$\gamma'^{\ast}=\gamma'^+-\gamma'^-$]{\includegraphics[scale=0.8]{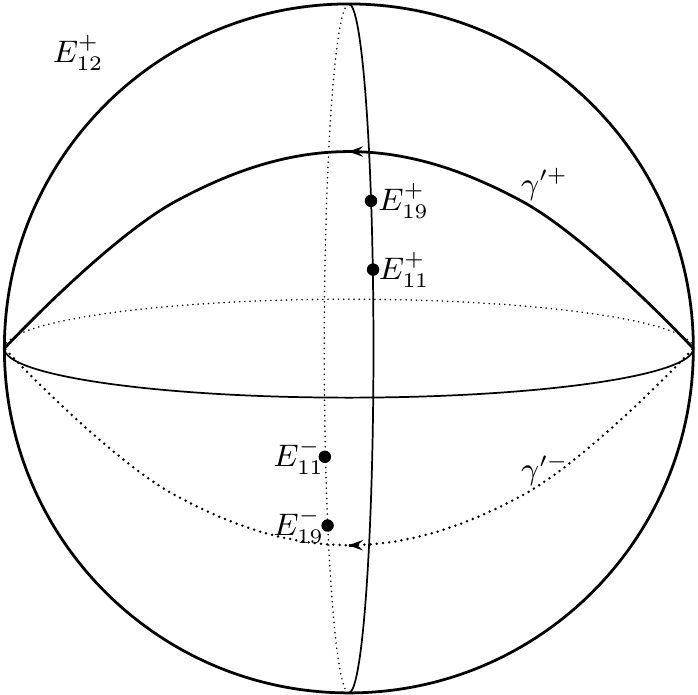}}
    \qquad
    \subfloat[$\gamma'^{\ast}=0$]{\includegraphics[scale=0.8]{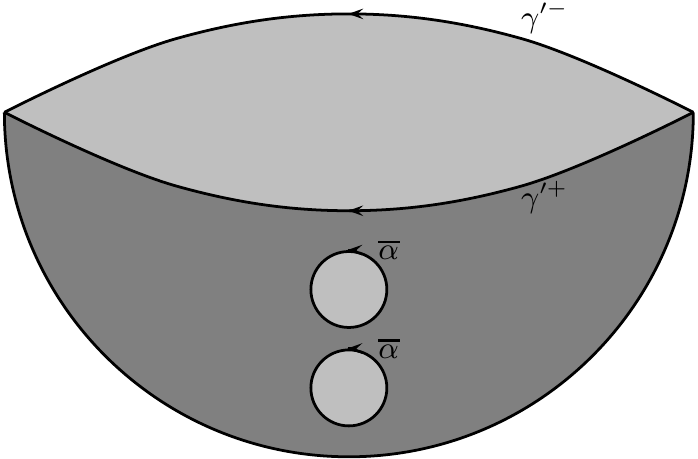}}
    \caption{$E_{12}^+$}
    \label{figure:E12+}
    \end{figure}

\begin{proof}
We construct a real $2$-dimensional surface $U$ in $\mathbb{P}^1 \times \mathbb{P}^1$ such that the path $\gamma$ is a boundary component of $U$. Define a path $\gamma'$ lying on $F_{1-}$ which connects the two points $F_{1-} \cap L_{2+}=\{(-1/\sqrt{-3}:1,1:1)\}$ and $F_{1-} \cap L_{2-}=\{(-1/\sqrt{-3}:1,-1:1)\}$ as follows: Let us define paths $\gamma_i'$ on $F_{1-}$ by
    \begin{align*}
    \gamma_1'(s) &= (-1/\sqrt{-3}:1, 1+s\sqrt{-1}:1), 0 \le s \le 2, \\
    \gamma_2'(s) &= (-1/\sqrt{-3}:1, -s+2\sqrt{-1}:1), -1 \le s \le 1, \\
    \gamma_3'(s) &= (-1/\sqrt{-3}:1, -1+(2-s)\sqrt{-1}: 1), 0 \le s \le 2,
    \end{align*}
and the path $\gamma'$ is defined by their union, that is,
    \begin{equation*}
    \gamma' = \gamma_1' \cup \gamma_2' \cup \gamma_3'.
    \end{equation*}
And then the surface $U$ is defined by connecting the paths $\gamma$ and $\gamma'$ by line segments. More precisely,
    \begin{equation*}
    \small
    \begin{split}
    U &= (t\gamma_1(s) + (1-t)\gamma_1'(s)) \cup (t\gamma_2(s) + (1-t)\gamma_2'(s)) \cup (t\gamma_2(s) + (1-t)\gamma_2'(s))\\
    &=\{(t(2+2s\sqrt{-1})+(1-t)(-1/\sqrt{-3}):1, 1+s\sqrt{-1}:1) \mid 0 \le s \le 2\} \\
    &\cup \{(t(-2s+4\sqrt{-1})+(1-t)(-1/\sqrt{-3}):1, -s+2\sqrt{-1}:1) \mid -1 \le s \le 1\} \\
    &\cup \{(t(-2+2(2-s)\sqrt{-1})+(1-t)(-1/\sqrt{-3}):1, -1+(2-s)\sqrt{-1}: 1) \mid 0 \le s \le 2\},
    \end{split}
    \end{equation*}
where $0 \le t \le 1$; Figure~\ref{figure:UandU'}(A).

    \begin{figure}[htb]
    \centering
    \subfloat[$U$]{\includegraphics[scale=0.9]{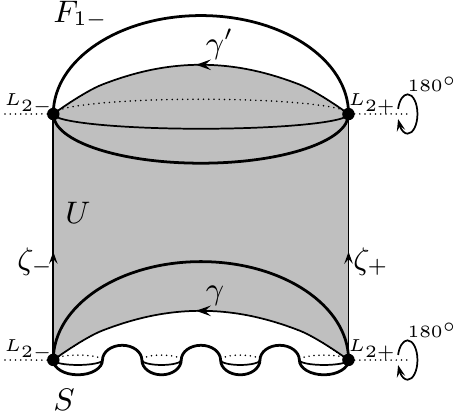}}
    \qquad
    \subfloat[$U_T$]{\includegraphics[scale=0.9]{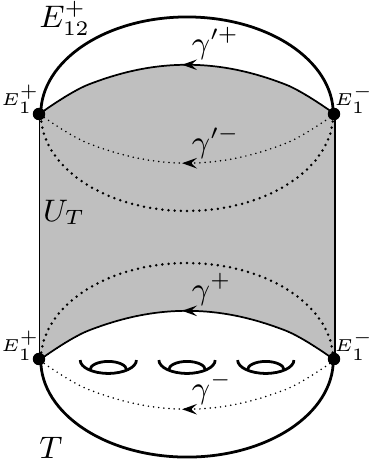}}
    \caption{The real $2$-dimensional surfaces $U$ and $U_T$}
    \label{figure:UandU'}
    \end{figure}

The surface $U$ does not pass through the blowing points of $R=(\mathbb{P}^1 \times \mathbb{P}^1) \sharp 16 \overline{\mathbb{P}^2} \to \mathbb{P}^1 \times \mathbb{P}^1$. Hence one may assume that $U$ is contained in $R$. The boundary $\partial U$ of $U$ consists of four components $\gamma$, $\gamma'$, $\zeta_-$, $\zeta_+$, where $\zeta_{\pm}$ is the line segment connecting two points $S \cap L_{2\pm}$ and $F_{1-} \cap L_{2\pm}$, respectively; that is,
    \[\zeta_{\pm} = \{(2t+(1-t)(-1/\sqrt{-3}):1, \pm 1:1) \mid 0 \le t \le 1\}.\]

Note that the curves $L_{2\pm}'$ are components of the branch divisors of the double covering $\phi : V \to R$. But the two boundary components $\zeta_-$ and $\zeta_+$ of $U$ are contained in $L_{2\pm}'$. Therefore the inverse image $U_T=\phi^{-1}(U)$ has only two boundaries:
    \[\partial U_T = \phi^{-1}(\gamma) \cup \phi^{-1}(\gamma').\]
Moreover, since the inverse image of the sphere $F_{1-}$ of the double covering $\phi$ is again the sphere $E_{12}^+$, it follows that the loop $\phi^{-1}(\gamma')=\gamma'^+-\gamma'^-$ lies in $E_{12}^+$; Figure~\ref{figure:UandU'}(B). In fact, by a simple calculation, one can show that the interior of $U$ does not intersect with the $(4,4)$-divisor $B=C_{\pm}+L_{1\pm}+L_{2\pm}$. Therefore $U_T$ is topologically a cylinder with two boundaries $\phi^{-1}(\gamma)$ and $\phi^{-1}(\gamma')$.

Set $\gamma'^{\ast} = \gamma'^+-\gamma'^-$. It follows by a simple calculation that $U \cap (L_{1-} \cup L_{1+} \cup G) = \varnothing$ and $U \cap L_3 = \{(2\sqrt{-1}:1, 2\sqrt{-1}:1)\}$. Therefore $U_T \cap (E_3^+ \cup E_3^-) = U_T \cap (E_{19}^+ \cup E_{19}^-) = \varnothing$ and $U_T \cap (E_{16}^{+} \cup E_{16}^{-})$ consists of two points. Since $U_T$ does not pass through any blowing-up points of $W=V \sharp 10 \overline{\mathbb{P}^2}$, $U_T$ does not intersect $\widetilde{E}_2^{\pm}$, $\widetilde{E}_4^{\pm}$, $\widetilde{E}_5^{\pm}$, $\widetilde{E}_6^{\pm}$, $\widetilde{E}_{11}^{\pm}$. Therefore $U_T' = U_T \cap W_0$ has four boundaries as in Figure~\ref{figure:Uitself}.
Since $2\overline{\alpha}=0$, we have
    \begin{equation}\label{equation:gamma+gamma'=0}
    \gamma^{\ast}+\gamma'^{\ast}=0.
    \end{equation}

    \begin{figure}[htb]
    \centering
    \includegraphics[scale=0.9]{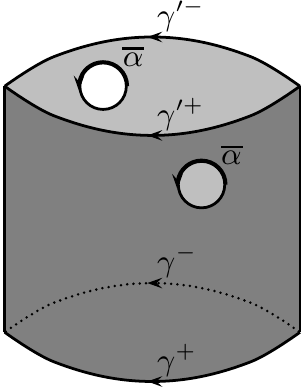}
    \caption{Four boundaries of $U_T' = U_T \cap W_0$}
    \label{figure:Uitself}
    \end{figure}

On the other hand the sphere $F_{1-}'$ intersects transversely with $H_3'$ at $(-1/\sqrt{-3}:1,-1/\sqrt{-3}:1)$ and with $G'$ at $(-1/\sqrt{-3}:1,\sqrt{-3}:1)$; Figure~\ref{figure:F1-}(A). Therefore the sphere $E_{12}^+$ intersects transversely with $E_{11}^+$ and $E_{19}^+$ on the upper half side of $E_{12}^+$ and with $E_{11}^-$ and $E_{19}^-$ on the lower half side of $E_{12}^+$; Figure~\ref{figure:E12+}(A). Therefore the loop $\gamma'^{\ast}$ encloses the two intersection points $E_{12}^+ \cap E_{11}^+$ and $E_{12}^+ \cap E_{19}^+$. The normal circles over the sphere bundles $C_2^+$ and $C_3^+$ over $E_{11}^+$ and $E_{19}^+$ are $\overline{\alpha}$. Therefore in $W_0$ the lower half sphere of $\widetilde{E}_{12}^+$ has three boundaries: $\gamma'^{\ast}$ and two $\overline{\alpha}$'s on the two intersection points $E_{12}^+ \cap E_{11}^+$ and $E_{12}^+ \cap E_{19}^+$; Figure~\ref{figure:E12+}(B). Since $2\overline{\alpha}=0$, we have $\gamma'^{\ast} = 0$. Therefore the relation~\eqref{equation:gamma+gamma'=0} implies that $\gamma^{\ast}=0$.
\end{proof}

\begin{lemma}\label{lemma:beta}
The loop $\beta$ is a generator of $H_1(Z_0;\mathbb{Z})$; hence, $H_1(Z_0;\mathbb{Z}) = \mathbb{Z}/4\mathbb{Z}$.
\end{lemma}

\begin{proof}
By the previous lemmas, we have $\widetilde{2\beta}=\overline{\alpha} \neq 0$ in $H_1(W_0;\mathbb{Z})$. Hence we get $2\beta = \alpha \neq 0$ in $H_1(Z_0;\mathbb{Z})$. But we have $H_1(Z_0;\mathbb{Z})=\langle \alpha, \beta \rangle$ by Lemma~\ref{lemma:<alpha,beta>}. Therefore $H_1(Z_0;\mathbb{Z})$ is generated by the single element $\beta$. Since $H_1(Z_0;\mathbb{Z})$ is of order $4$ by Lemma~\ref{lemma:|H_1(Z_0)|=4}, we have $H_1(Z_0;\mathbb{Z})=\mathbb{Z}/4\mathbb{Z}$.
\end{proof}

We finally prove the main result of this section.

\begin{proposition}\label{proposition:H_1(Zbar)=Z_4}
$H_1(\overline{Z};\mathbb{Z})=\mathbb{Z}/4\mathbb{Z}$.
\end{proposition}

\begin{proof}
We consider the Mayer-Vietoris sequence of a pair
$(Z_0, M_{1} \cup M_{2} \cup  M_{3})$:
    \begin{equation*}\label{eqation:MV(Z,M)}
    H_1(\partial M_{1} \cup \partial M_{2} \cup \partial M_{3};\mathbb{Z})
    \xrightarrow{k_{\ast} \oplus j_{\ast}} H_1(Z_0;\mathbb{Z}) \oplus H_1(M_{1} \cup M_{2} \cup M_{3};\mathbb{Z}) \to H_1(\overline{Z};\mathbb{Z}) \to 0.
    \end{equation*}
Let $(1,0,0), (0,1,0)$ and $(0,0,1)$ be generators of
    \begin{equation*}
    \begin{split}
    H_1(\partial M_{1} \cup \partial M_{2} \cup \partial M_{3};\mathbb{Z}) & \cong H_1(\partial M_{1};\mathbb{Z})\oplus H_1(\partial M_{2};\mathbb{Z}) \oplus H_1(\partial M_{3};\mathbb{Z}) \\
    & \cong \mathbb{Z}/72 \mathbb{Z} \oplus \mathbb{Z}/4 \mathbb{Z} \oplus \mathbb{Z}/4 \mathbb{Z}.
    \end{split}
    \end{equation*}
Then the images of $(1,0,0), (0,1,0)$ and $(0,0,1)$ under the map $k_{\ast} \oplus j_{\ast}$ are $(\alpha,1,0,0)$, $(\alpha, 0,1,0)$ and $(\alpha,0,0,1)$, respectively, in
    \begin{equation*}
    \begin{split}
    &H_1(Z_0;\mathbb{Z}) \oplus H_1(M_{1} \cup M_{2} \cup M_{3};\mathbb{Z})\\
    & \qquad \qquad \qquad \qquad \cong H_1(Z_0;\mathbb{Z}) \oplus H_1(M_{1};\mathbb{Z}) \oplus H_1(M_{2};\mathbb{Z}) \oplus H_1(M_{3};\mathbb{Z})\\
    & \qquad \qquad \qquad \qquad \cong \mathbb{Z}/4\mathbb{Z} \oplus \mathbb{Z}/6\mathbb{Z} \oplus \mathbb{Z}/2\mathbb{Z}\oplus \mathbb{Z}/2\mathbb{Z}.
    \end{split}
    \end{equation*}
Since $\alpha=2\beta$ and $\beta$ is a generator of $H_1(Z_0;\mathbb{Z}) \cong \mathbb{Z}/4\mathbb{Z}$, it follows that
    \begin{equation*}
    \begin{split}
    H_1(\overline{Z}; \mathbb{Z}) & \cong H_1(Z_0;\mathbb{Z}) \oplus H_1(M_{1} \cup M_{2}
    \cup M_{3};\mathbb{Z}) /\im(k_{\ast} \oplus j_{\ast})  \\
    & \cong (\mathbb{Z}/4\mathbb{Z} \oplus \mathbb{Z}/6\mathbb{Z} \oplus \mathbb{Z}/2\mathbb{Z}\oplus \mathbb{Z}/2\mathbb{Z})/ \langle (2,1,0,0),(2,0,1,0),(2,0,0,1)\rangle \\
    & \cong \langle (1,0,0,0)\rangle \cong\mathbb{Z}/4\mathbb{Z}.
    \end{split}
    \end{equation*}
Therefore we finally get $H_1(\overline{Z};\mathbb{Z}) \cong \mathbb{Z}/4\mathbb{Z}$.
\end{proof}

\section{More examples}\label{section:more_example}

In this section we construct three more examples.

\subsection{An example with an ample canonical divisor}\label{subsection:main-modified-ample}

We construct an example with an ample canonical divisor by modifying the configuration in the main construction, while the main example in Section~\ref{section:main-construction} may have a non-ample canonical divisor as stated in Proposition~\ref{proposition:main-not-ample}.

In order to construct a singular surface $X$ with three singularities of class $T$ in Section~\ref{section:main-construction}, we contracted three linear chains $C_1$, $C_2$, $C_3$ of rational curves from the blown-up surface $Z$ so that we obtained a singular surface $X$ with three singularities of class $T$; cf.\thinspace Figure~\ref{figure:Z}. Together with $C_i$'s, contract one more rational curve, the proper transform of $E_{10}$. Since the proper transform of $E_{10}$ is a ($-2$)-curve, we obtain a projective singular surface $X'$ with four singularities of class $T$. Let $f: Z \to X'$ be the contraction. By similar methods in the previous Section~\ref{section:global-Q} and Section~\ref{section:Z_4} one can show that there is a global a $\mathbb{Q}$-Gorenstein smoothing of $X'$ and its general fiber $X_t'$ is a minimal surface of general type with $p_g=0$, $K^2=2$, and $H_1=\mathbb{Z}/4\mathbb{Z}$.

\begin{proposition}\label{proposition:main-modified-ample}
A general fiber $X_t'$ of a $\mathbb{Q}$-Gorenstein smoothing of $X'$ has an ample canonical divisor.
\end{proposition}

\begin{proof}
We use a similar method in the author's paper~\cite{PPS-K4} and J. Keum-Y. Lee-H. Park~\cite{Keum-Lee-Park}. We denote by
    \[C_1: \udc{-7}{u_1}-\udc{-3}{u_2}-\udc{-2}{u_3}-\udc{-2}{u_4}-\udc{-2}{u_5}-\udc{-2}{u_6}, \  C_2: \udc{-4}{u_7}, \  C_3: \udc{-4}{u_8}, \ C_4: \udc{-2}{u_9}\]
the linear chains of rational curves, where $u_1=\widetilde{E}_{16}$,  $u_2=\widetilde{E}_2$, $u_3=\widetilde{E}_3$, $u_4=\widetilde{E}_4$, $u_5=\widetilde{E}_5$, $u_6=\widetilde{E}_6$, $u_7=\widetilde{E}_{11}$, $u_8=\widetilde{E}_{19}$, $u_9=\widetilde{E}_{10}$ denote the proper transforms of $E_i$ on the blown-up space $Z$. By a similar calculation in the proof of Theorem~\ref{theorem:properties_of_X_t}, it follows that
    \begin{equation}\label{equation:f^*K_X'}
    \begin{split}
    f^{\ast}{K_{X'}} \equiv &\frac{5}{6} u_1 + \frac{5}{6} u_2 + \frac{4}{6} u_3 + \frac{3}{6} u_4 + \frac{2}{6} u_5 + \frac{1}{6} u_6 + \frac{1}{2} u_7 + \frac{1}{2} u_8\\
    &+ e_1 + e_2 + e_3 + e_4 + e_5
    \end{split}
    \end{equation}
and that $f^{\ast}{K_{X'}}$ is nef and, furthermore, $f^{\ast}{K_{X'}} \cdot e_j > 0$ for all $j=1,\dotsc,5$.

Since $K_{X'}^2=2>0$, in order to show the ampleness of $K_{X'}$, it is enough to show that $K_{X'} \cdot C > 0$ for any irreducible curve $C \subset X'$. Let $\overline{C} \subset Z$ be the proper transform of $C$. Note that
    \[K_{X'} \cdot C = f^{\ast}{K_{X'}} \cdot f^{\ast}{C} = f^{\ast}{K_{X'}} \cdot \overline{C}.\]
Hence it is enough to show that $f^{\ast}{K_{X'}} \cdot \overline{C} > 0$.
But we already know that $f^{\ast}{K_{X'}} \cdot e_j > 0$ for all $j=1,\dotsc,5$. Thus we may assume that $\overline{C} \neq e_j$ for any $j$. Since the coefficients of $f^{\ast}{K_{X'}}$ in \eqref{equation:f^*K_X'} are positive and $\overline{C} \neq u_i$ for any $i=1,\dotsc,9$ and $\overline{C} \neq e_j$ for any $j=1,\dotsc,5$, if $\overline{C} \cdot u_i > 0$ for some $i$ or $\overline{C} \cdot e_j > 0$ for some $j$, then $f^{\ast}{K_{X'}} \cdot \overline{C} > 0$.

Let $h: Z \to Y$ be the blowing-up that produces the exceptional curve $e_j$'s. We denote $h(\overline{C})$ by $\widehat{C}$. If $p_a(\overline{C}) \ge 2=p_a(\widehat{C}) \ge 2$, then $\widehat{C} \cdot (E_{16}+E_{19}) > 0$. Therefore $\overline{C} \cdot h^{\ast}(E_{16}+E_{19}) > 0$. If $p_a(\overline{C})=p_a(\widehat{C})=1$, then $\widehat{C}$ is a fiber or a half pencil of an elliptic pencil. If $\widehat{C}$ is numerically equivalent to $\lvert E_{16} + E_{19} \rvert$ or $(1/2)(E_{16}+E_{19})$, then we have $\widehat{C} \cdot E_2 > 0$ because $E_2$ is a bisection; hence, $\overline{C} \cdot h^{\ast}{E_2} > 0$. If not, then $\widehat{C} \cdot (E_{16} + E_{19}) > 0$. Finally, if $p_a(\overline{C})=0$ then $\overline{C}$ is a rational curve, hence $\overline{C} = \widetilde{E}_k$ for some $1 \le k \le 20$. However, according to the dual graph (Figure~\ref{figure:dual_graph}) of (all) rational curves lying on the Enriques surface $Y$, every rational curve lying on $Z$ which is not contracted to a singular points of $X'$ intersects with some $u_i$'s. In any case we show that $\overline{C} \cdot u_i > 0$ for some $i$ or $\overline{C} \cdot e_j > 0$ for some $j$, which completes the proof.
\end{proof}

\subsection{An example obtained from a different configuration}\label{subsection:more-not-ample}

We construct two more examples of a minimal surface of general type with $p_g=0$, $K^2=2$ and $H_1=\mathbb{Z}/4\mathbb{Z}$ using different configurations of rational curves coming from the same Enriques surface $Y$. We first construct an example with a non-ample canonical divisor; Remark~\ref{remark:more-example-non-ample}. We then construct another example by modifying the configuration. But we don't know whether the canonical divisor of the latter example is ample or not; Remark~\ref{remark:more-example-ample-or-not}. Since all proofs are basically the same as the case of the main example constructed in Section~\ref{section:main-construction}, we only explain how to construct the blown-up surface $Z$ from the Enriques surface $Y$ containing linear chains of rational curves.

\begin{figure}[htb]
    \centering
    \subfloat[An Enriques surface $Y$]{\includegraphics[scale=0.9]{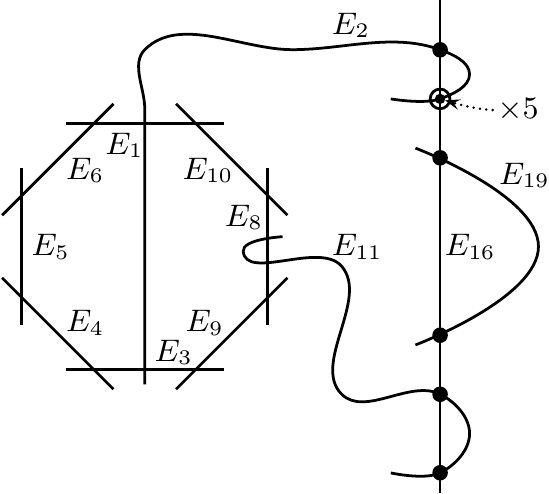}}
    \qquad
    \subfloat[$Z=Y \sharp 10 \overline{\mathbb{P}^2}$]{\includegraphics[scale=0.9]{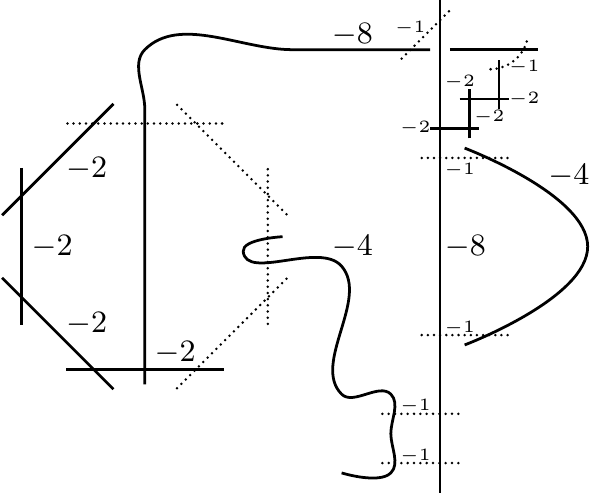}}
    \caption{Another example with $p_g=0$, $K^2=2$ and $H_1=\mathbb{Z}/4\mathbb{Z}$}
    \label{figure:another-example}
    \end{figure}

On the Enriques surface $Y$ used in Section~\ref{section:main-construction}
we blow up 5 times at the marked points $\bullet$ on $E_{16}$ and we blow up again successively $5$ times at the marked point $\bigodot$ on $E_{16}$ ; Figure~\ref{figure:another-example}(A). The blown-up surface $Z=Y \sharp 10 \overline{\mathbb{P}^2}$ contains the following four disjoint linear chains of rational curves as in Figure~\ref{figure:another-example}(B):
    \begin{align*}
    &C_1=C_{6,1}: \udc{-8}{u_1}-\udc{-2}{u_2}-\udc{-2}{u_3}-\udc{-2}{u_4}-\udc{-2}{u_5}, \\
    &C_2=C_{6,1}: \udc{-8}{u_6}-\udc{-2}{u_7}-\udc{-2}{u_8}-\udc{-2}{u_9}-\udc{-2}{u_{10}},\\
    &C_3=C_{2,1}: \udc{-4}{u_{11}}, \quad C_4=C_{2,1}: \udc{-4}{u_{12}},
    \end{align*}
where $C_1$ consists of the proper transforms of $E_2$, $E_3$, $E_4$, $E_5$, $E_6$, and $C_2$, $C_3$, $C_4$ contain the proper transforms of $E_{16}$, $E_{11}$, $E_{19}$, respectively.

We contract the four linear chains $C_1$, $C_2$, $C_3$, $C_4$ from $Z$. Then we obtain a projective singular surface $X$ with four singularities of class $T$. By applying the same method to the singular surface $X$ as in the previous sections, one can show that a general fiber $X_t$ of a $\mathbb{Q}$-Gorenstein smoothing of $X$ is a minimal surface of general type with $p_g=0$, $K^2=2$ and $H_1=\mathbb{Z}/4\mathbb{Z}$.

\begin{remark}\label{remark:more-example-non-ample}
By applying the same method in the proof of Proposition~\ref{proposition:main-not-ample}, we can show that there is a $\mathbb{Q}$-Gorenstein smoothing of the singular surface $X$ such that the canonical divisor $K_{X_t}$ of a general fiber $X_t$ of the $\mathbb{Q}$-Gorenstein smoothing of $X$ is \emph{not} ample.
\end{remark}

\begin{remark}\label{remark:more-example-ample-or-not}
As in the previous subsection, let us contract one more $(-2)$-curve $\widetilde{E}_{10}$ as well as four linear chains $C_1$, $C_2$, $C_3$, $C_4$ from $Z$. Let $f: Z \to X'$ be the contraction. One can show that there is a global a $\mathbb{Q}$-Gorenstein smoothing of $X'$ and its general fiber $X_t'$ is a minimal surface of general type with $p_g=0$, $K^2=2$ and $H_1=\mathbb{Z}/4\mathbb{Z}$.
However we cannot conclude that $K_{X'}$ is ample because, letting $e$ be the exceptional divisor connecting $u_1$ and $u_6$, we have $f^{\ast}{K_{X'}} \cdot e = 0$. We don't know whether the canonical divisor $K_{X_t'}$ of a general fiber $X_t'$ is ample or not. We leave this problem for future research.
\end{remark}

\section{Appendix: An example with $H_1=\mathbb{Z}/2\mathbb{Z} \oplus \mathbb{Z}/2\mathbb{Z}$}\label{section:appendix}

As mentioned in the introduction, minimal complex surfaces with $p_g=0$, $K^2=2$ and $H_1=\mathbb{Z}/2\mathbb{Z} \oplus \mathbb{Z}/2\mathbb{Z}$ have been constructed in Inoue~\cite{Inoue} and Keum~\cite{Keum} by the classical method: Quotient by group action. In this appendix we construct such a surface by a $\mathbb{Q}$-Gorenstein smoothing method. We explain briefly how to construct such an example and we sketch how to prove $H_1=\mathbb{Z}/2\mathbb{Z} \oplus \mathbb{Z}/2\mathbb{Z}$ because almost all proofs are basically the same as the case of $H_1=\mathbb{Z}/4\mathbb{Z}$.

\subsection*{Construction}

\begin{figure}[htb]
    \centering
    \subfloat[An Enriques surface $Y$]{\includegraphics[scale=0.9]{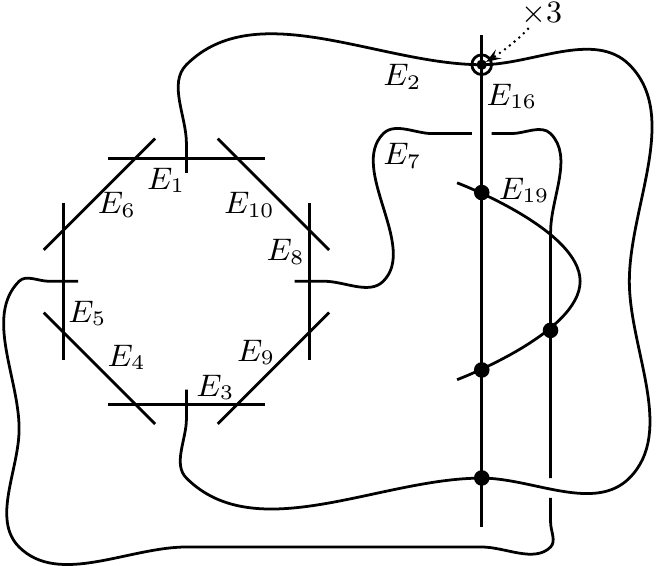}} \quad
    \subfloat[$Z=Y \sharp 7 \overline{\mathbb{P}^2}$]{\includegraphics[scale=0.9]{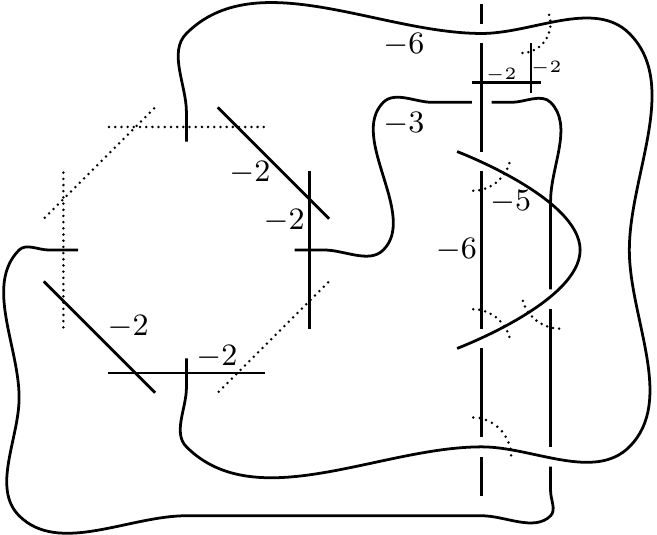}}
    \caption{An example with $p_g=0$, $K^2=2$ and $H_1=\mathbb{Z}/2\mathbb{Z} \oplus \mathbb{Z}/2\mathbb{Z}$}
    \label{figure:Z2+Z2}
    \end{figure}

We begin with the same Enriques surface $Y$ used in Section~\ref{section:main-construction} for constructing an example with $H_1=\mathbb{Z}/4\mathbb{Z}$. However we use another bisection $E_{7}$ instead of the bisection $E_{11}$; Figure~\ref{figure:Z2+Z2}(A). We blow up the four marked point $\bullet$ and we blow up again $3$ times at the marked point $\bigodot$. The blown-up surface $Z=Y \sharp 7 \overline{\mathbb{P}^2}$ has four disjoint linear chains of rational curves as in Figure~\ref{figure:Z2+Z2}(B):
    \begin{equation*}
    C_1:\uc{-5}-\uc{-3}-\uc{-2}-\uc{-2}, \ \ \ C_2: \uc{-6}-\uc{-2}-\uc{-2},
    \ \ \ C_3: \uc{-6}-\uc{-2}-\uc{-2},
    \end{equation*}
where $C_1$ consists of the proper transforms of $E_{19}$, $E_7$, $E_8$, $E_{10}$ and $C_2$ consists of the proper transforms of $E_2$, $E_3$, $E_4$. The chain $C_3$ contains the proper transforms of $E_{16}$.

By contracting these three disjoint chains from $Z$, we get a projective surface $X$ with three singularities of class $T$. Then, by the same argument in Section~\ref{section:main-construction} and Section~\ref{section:global-Q}, we see that the surface $X$ has a $\mathbb{Q}$-Gorenstein smoothing and a general fiber $X_t$ of the $\mathbb{Q}$-Gorenstein smoothing for $X$ is a minimal complex surface of general type with $p_g=0$ and $K^2=2$.

\begin{proposition}
$H_1(X_t, \mathbb{Z}) = \mathbb{Z}/2\mathbb{Z} \oplus \mathbb{Z}/2\mathbb{Z}$.
\end{proposition}

\begin{proof}
Let $\overline{Z}$ be the rational blow-down $4$-manifold obtained from $Z$ by replacing $C_i$'s with the corresponding Milnor fibers, respectively. As in the case of the proof of $H_1=\mathbb{Z}/4\mathbb{Z}$, it is enough to show that $H_1(\overline{Z};\mathbb{Z})=\mathbb{Z}/2\mathbb{Z} \oplus \mathbb{Z}/2\mathbb{Z}$. We decompose $Z = Z_0 \cup \{C_1 \cup C_2 \cup C_3\}$. Applying the same proof of Proposition~\ref{proposition:H_1(Zbar)=Z_4} to this case, one can show that $H_1(\overline{Z};\mathbb{Z})=\mathbb{Z}/2\mathbb{Z} \oplus \mathbb{Z}/2\mathbb{Z}$ if $H_1(Z_0;\mathbb{Z})=\mathbb{Z}/2\mathbb{Z} \oplus \mathbb{Z}/2\mathbb{Z}$.

We consider the loop $\beta$ on the curve $T/\sigma$ constructed in \S\ref{subsubsection:beta}. In Lemma~\ref{lemma:beta} we proved that if $\widetilde{2\beta} \neq 0$ then  $H_1(Z_0;\mathbb{Z})=\mathbb{Z}/4\mathbb{Z}$. By a similar argument, one can show that if $\widetilde{2\beta} = 0$ then $H_1(Z_0;\mathbb{Z})=\mathbb{Z}/2\mathbb{Z} \oplus \mathbb{Z}/2\mathbb{Z}$.

In order to prove $\widetilde{2\beta} = 0$, we consider the real $2$-dimensional surface $U_T'$ constructed in Lemma~\ref{lemma:gamma^*=0}. Since the proper transform of $E_{16}$ is removed to obtain $Z_0$, one may apply the same proof of Equation~\eqref{equation:gamma+gamma'=0} to our case, so that we get
    \[\gamma^{\ast} + \gamma'^{\ast} = 0.\]
On the other hand we use $E_7$ instead of $E_{11}$ to construct the linear chains $C_i$. Since $E_7$ does not intersects with $E_{12}$, we have only one $\alpha$ and $\gamma'^{\ast}$ as boundary in Figure~\ref{figure:E12+}(B). Hence we have
    \[\gamma'^{\ast}+\overline{\alpha}=0.\]
Since $\widetilde{2\beta}=\gamma^*+\overline{\alpha}$ by Lemma~\ref{lemma:2beta=gamma^*+alpha}, combining these relations, we have $\widetilde{2\beta}=0$.
\end{proof}

\providecommand{\bysame}{\leavevmode\hbox to3em{\hrulefill}\thinspace}

\end{document}